\documentclass[11pt]{amsart}
\usepackage[francais]{babel}
\usepackage{t1enc,amssymb}
\usepackage[author-year]{amsrefs}
\newcommand\R{{\mathbb R}}
\newcommand\fcite[3]{\textsf{#1\footnote{\cite{#2}*{page #3}}}}

\begin{document}
\title{Théorie des groupes et psychologie de l'intelligence}
\author{Laurent Bartholdi}
\author{Fabrice Liardet}
\date{21 février 2005; composé \today}
\maketitle
\tableofcontents

\section{Introduction}
Le but de cette note est de rapporter des résultats sur l'intersection entre la
théorie des groupes et la psychologie de l'intelligence, telle qu'elle a été
étudiée par Piaget et son école.

Piaget a fait un usage étendu d'idées et de techniques mathématiques,
la plus frappante étant peut-être son utilisation de la théorie des
catégories~\cite{piaget:mc}.  Toutefois, comme S. Pappert nous met en
garde dans la préface de cet ouvrage, ``\fcite{Ma propre expérience de
  ce livre n'a pas été sans quelque douleur\dots puis-je voir dans
  l'utilisation que fait Piaget d'idées mathématiques tout à fait
  profondes autre chose qu'une métaphore superficielle?  Je ne suis
  pas le seul à être parfois assailli de tels doutes.}{piaget:mc}{8}''
La source du problème est facile à comprendre: Piaget écrit pour des
psychologues, dans un langage que les mathématiciens ont beaucoup de
peine à comprendre; et il utilise des concepts naturels pour les
mathématiciens, mais en court-circuitant les années d'étude
nécessaires à un étudiant de mathématiques pour se familiariser à ces
objets.

Un projet ambitieux serait de réécrire les ouvrages de Piaget (40 volumes) de
façon qu'ils soient plus accessibles aux mathématiciens; le vocabulaire
mathématique de Piaget est en effet non-standard, et le contenu psychologique y
est très dense. Nous nous proposons, comme but plus modeste, de recenser les
usages de théorie des groupes dans la théorie piagetienne, et ainsi de fournir
à la fois une introduction de Piaget aux mathématiciens et une introduction des
mathématiques aux piagetiens.

La théorie des groupes semble en effet le domaine \emph{par
  excellence} des mathématiques propre à décrire le monde ---
ainsi, par exemple, les déplacements sont régis par un groupe,
comme on le verra au \S\ref{ss:gd}. L'axiome le plus important de la
théorie des groupes, l'\emph{associativité}, contient en
lui-même l'idée d'abstraction. Nous verrons donc en passant les
divers domaines non-mathématiques où la théorie des groupes
s'applique.

Finalement, à titre d'illustration, et dans le pur esprit piagetien, nous
décrirons quelques jeux, dont certains ont déjà été implémentés sur ordinateur,
qui font un usage non-trivial de la théorie des groupes. Un soin particulier a
été pris pour qu'aucun concept mathématique ne soit apparent dans ces jeux.
Notre espoir est qu'un enfant, y jouant, découvre les \emph{idées}
des mathématiques sous-jacentes, sans y voir la forme.

\subsection{Remerciements}
Nous tenons à exprimer notre immense gratitude à Androula
Henriques-Christophidès et à Roger Flückiger, qui par leurs
discussions généreuses nous ont grandement facilité l'accès aux textes
de Piaget. Nous sommes aussi reconnaissants à Daniela Helbig, Vaughan
Jones et Florian Sobiecky de nous avoir expliqué certaines relations
entre la physique théorique et la théorie des groupes.

\section{Introduction à la théorie des groupes}\label{s:intromath}
Le lecteur trouvera ci-dessous un bref résumé des notions importantes en
théorie des groupes.

Un groupe peut être défini de deux manières différentes, comme un
groupe \emph{abstrait} ou \emph{concret}; un résultat fondamental est
que ces deux définitions sont équivalentes.  Commençons par les
groupes abstraits; dans ces groupes, les éléments sont représentés par
des symboles.  Soit $X$ un ensemble. On dira que $X$ est muni d'une
structure de groupe s'il existe (1) une opération, nommée
\emph{multiplication}\footnote{qui n'a en général rien à voir avec la
  multiplication habituelle des entiers.}, qui à deux éléments $x,y\in
X$ attribue un produit $x\cdot y\in X$, (2) une opération, nommée
\emph{inversion}, qui à $x\in X$ attribue un inverse $x^{-1}\in X$, et
(3) un \emph{élément neutre} $1\in X$, de sorte pour $x,y,z$ des éléments 
quelconques de $X$ les relations suivantes soient vérifiées:
\begin{alignat}{2}
  (xy)z &= x(yz) & \text{associativité}\label{eq:assoc}\\
  xx^{-1} &= x^{-1}x = 1 & \text{inverse}\label{eq:inv}\\
  x1 &= 1x = x & \text{neutre}\label{eq:id}
\end{alignat}
Ces trois propriétés sont les \emph{axiomes} définissant la structure
de groupe. Le premier axiome stipule que la \emph{multiplication} de
trois éléments peut se faire dans l'ordre qu'on veut; c'est-à-dire
qu'étant donnés trois éléments $x,y,z$ dans cet ordre, leur produit
$xyz$ peut se calculer en multipliant d'abord $x$ avec $y$, puis ce
résultat intermédiaire avec $z$, ou bien en multipliant d'abord $y$
avec $z$, puis en multipliant $x$ avec ce résultat intermédiaire.

Remarquons que l'associativité ne signifie \emph{pas} qu'on peut
mettre les trois éléments dans l'ordre qu'on veut quand on les
multiplie. Il existe une propriété additionnelle qu'un groupe peut
avoir, et qui dit que pour multiplier deux éléments $x,y$ on peut
calculer $xy$ ou $yx$, obtenant le même résultat dans les deux cas.
Cette propriété indique qu'on peut échanger l'ordre dans un produit de
deux (et donc un nombre quelconque d') élements, et s'appelle
\emph{commutativité}\footnote{L'associativité peut être vue comme une
  forme de commutativité à un niveau supérieur d'abstraction: elle
  indique que, étant donné un élément $y$, les opérations <<le
  multiplier par $x$ sur la gauche>> et <<le multiplier par $z$ sur la
  droite>> commutent.}; elle s'exprime par l'axiome $xy=yx$.

Le deuxième axiome stipule que tout élément $x$ a un inverse,
c'est-à-dire un élément noté $x^{-1}$ tel que le produit de $x$ par
son inverse donne l'élément neutre.  Le troisième axiome stipule qu'un
élément est inchangé si on le multiplie par l'élément neutre.

À titre d'exemple de groupe, considérons un ensemble $X$ à trois
éléments, qu'on notera $\cdot,\curvearrowright,\curvearrowleft$
respectivement. Alors une structure de groupe sur $X$ n'est rien
d'autre que la donnée de la multiplication, de l'inverse et du
neutre, sous forme de table par exemple:\footnote{On remarque que,
  dans la première définition de groupe, l'inverse et
  l'élément $1$ sont déterminés par la table de multiplication ---
  i.e., étant donné un groupe $X$ dont on ne connaît que la
  table de multiplication, on peut reconstruire l'opération
  $x^{-1}$ et l'élément $1$, et ce d'unique façon.}
\[\begin{array}{c|ccc}
  xy & y=\cdot & y=\curvearrowright & y=\curvearrowleft\\ \hline
  x=\cdot & \cdot & \curvearrowright & \curvearrowleft\\
  x=\curvearrowright & \curvearrowright & \curvearrowleft & \cdot\\
  x=\curvearrowleft & \curvearrowleft & \cdot & \curvearrowright
\end{array}\qquad
\begin{array}{c|ccc}
  x & \cdot & \curvearrowright & \curvearrowleft\\ \hline
  x^{-1} & \cdot & \curvearrowleft & \curvearrowright
\end{array}
\qquad\text{ élément neutre: }\cdot
\]
On remarque qu'on aurait aussi pu écrire
$\cdot=0,\curvearrowright=1,\curvearrowleft=2$, auquel cas la
multiplication s'interpréterait comme l'addition d'entiers modulo
$3$.

Un groupe concret est un groupe de transformations, défini comme suit.
Soit $U$ un ensemble. Un \emph{groupe sur $U$} est un ensemble $G$,
dont les éléments sont des bijections\footnote{Une bijection de $U$
  est une correspondance entre $U$ et $U$ qui relie chaque élément de
  $U$ à un unique élément de $U$.} de $U$, tel que si on prend deux
éléments quelconques $f,g\in G$ alors la fonction composée $f\circ g$
et la fonction inverse $f^{-1}$ appartiennent à $G$, et telle que la
transformation identique appartient à $G$ (cette dernière condition
est nécessairement satisfaite dès que $G$ est non-vide).

On appelle $G$ un groupe concret, car il contient de façon
sous-jacente des bijections d'un ensemble concret $U$. C'est pourtant
un objet abstrait, dans le sens où les éléments de $G$, qui sont
fondamentalement des bijections de $U$, sont considérés comme des
objets atomiques.

Reprenons l'exemple ci-dessus. Si on considère un triangle équilatéral
$U$ dans le plan, on peut le faire tourner d'un nombre quelconque de
tiers de tour pour obtenir à nouveau un triangle dans la même
position.  Cette rotation est une bijection du plan, et par
restriction une bijection des points du triangle, i.e.\ de $U$. La
composition et l'inverse de telles rotations est encore une rotation
du triangle. Si on tourne le triangle d'un tour entier, on n'a déplacé
aucun des points du triangle; ainsi par exemple la rotation d'un tiers
de tours est égale à la rotation de quatre tiers de tour. Il n'y a que
trois éléments dans ce groupe de rotations $G$. Si on les appelle
$\cdot,\curvearrowright,\curvearrowleft$ pour désigner respectivement
la transformation identique, la rotation d'un tiers de tour dans le
sens des aiguilles d'une montre, et la rotation d'un tiers de tour en
sens inverse, on voit que ces éléments se multiplient et s'inversent
comme décrit plus haut.

De façon plus générale, on peut toujours
\emph{représenter} un groupe abstrait $X$ comme un ensemble $G$ de
bijections de $U$. D'autre part, un ensemble $G$ de bijections de
$U$ satisfait toujours les axiomes définissant un groupe abstrait.
Ces définitions sont donc équivalentes --- ce résultat est connu
sous le nom de \emph{théorème de Cayley}.

\subsection{Catégories}\label{ss:cat} Dans son étude des structures
mentales construites pendant le développement de l'enfant, Piaget
considère aussi diverses théories plus faibles (i.e.\ satisfaisant des
axiomes plus faibles que ceux des groupes). Il ne les nomme pas
explicitement, bien qu'un riche vocabulaire mathématique existe. Un
\emph{semigroupe} est un ensemble $X$ muni d'une opération binaire et
satisfaisant~\eqref{eq:assoc}. Un \emph{monoïde} est un semigroupe
muni d'un élément neutre $1$ et satisfaisant~\eqref{eq:assoc}.  Un
\emph{magma} est un ensemble muni d'une opération binaire; ainsi
groupes, semigroupes et monoïdes sont des magmas.

On peut aussi relaxer la définition de la multiplication, et arriver ainsi à la
notion de catégorie et de groupoïde. Une \emph{catégorie}\footnote{Pour les
  mathématiciens: on ne parlera ici que de petite catégorie.} est une paire
d'ensembles $X,G$ munie d'applications $s,t:G\to X$. On doit penser à $X$ comme
un ensemble de points, et à $G$ comme un ensemble de flèches reliant
un point à un autre; alors $s(g)$ et
$t(g)$ sont respectivement la source et le but de le la flèche $g$. On a une
multiplication partiellement définie sur $G$: si $g,h\in G$ sont tels que
$s(g)=t(h)$, alors ils ont un produit $gh\in G$; on a $s(gh)=s(h)$ et
$t(gh)=t(h)$. On a aussi des unités $1_x\in G$, telles que $s(1_x)=t(1_x)=x$,
et les axiomes
\begin{alignat}{2}
  (gh)k &= g(hk) & \text{associativité}\label{eq:cassoc}\\
  g1_{s(g)} &= 1_{t(g)}g = g & \text{identité}\label{eq:cid}\\
\end{alignat}
Si de plus tout $g\in G$ a un inverse, i.e.\ un élément $g^{-1}\in G$ avec
$s(g^{-1})=t(g)$ et $t(g^{-1})=s(g)$ et $gg^{-1}=1_{t(g)}$ et
$g^{-1}g=1_{s(g)}$, on appelle $G$ un \emph{groupoïde}.

On a, en résumé: une catégorie pour laquelle $X$ n'a qu'un élément est un
semigroupe. Une catégorie pour laquelle tout élément est inversible est un
groupoïde. Un groupoïde pour lequel $X$ n'a qu'un élément est un groupe. Un 
semigroupe pour lequel tout élément est inversible est un groupe.

Il existe une autre structure mathématique, celle de \emph{relation
  d'équivalence}. C'est une relation sur $X$, i.e.\ une propriété $x\sim y$,
qui est transitive ($x\sim y$ et $y\sim z$ implique $x\sim z$), symétrique
($x\sim y$ implique $y\sim x$) et réflexive ($x\sim x$). Ces axiomes présentent
une similarité formelle avec la structure de groupe. Ceci est rendu plus
explicite en constatant qu'une relation d'équivalence peut être décrite à
l'aide d'un groupoïde: on prend pour $G$ l'ensemble des couples $(x,y)\in
X\times X$ tels que $x\sim y$. Les fonctions $s,t$ sont respectivement définies
comme la projection sur la première et la seconde coordonnée de $G$. Le produit
est donné par $(y,z)(x,y)=(x,z)$ et l'inverse est donné par $(x,y)^{-1}=(y,x)$.
En ce sens, un groupoïde est une généralisation naturelle des groupes et des
  relations d'équivalence.

\section{La théorie des groupes chez Piaget}\label{s:gpespiaget}
Il n'est pas à priori évident que les structures mathématiques
décrites dans le chapitre~\ref{s:intromath} s'appliquent à l'étude du
développement de l'enfant. C'est justement une des intuitions
brillantes de Piaget que d'avoir perçu ce lien: \fcite{\dots si
  tardive qu'ait été la découverte de la notion de <<groupe>> en tant
  qu'être mathématique, ses propriétés les plus générales expriment en
  réalité certains des mécanismes les plus caractéristiques de notre
  intelligence.}{piaget:eg1}{45}

La théorie de Piaget s'appuie sur un grand nombre d'expériences réalisées sur
des enfants. Une des plus célèbres est la suivante: \fcite{Deux petits verres
  $A$ et $A_2$ de forme et de dimensions égales sont remplis d'un même nombre
  de perles, cette équivalence étant reconnue par l'enfant qui les a lui-même
  placées, par exemple en mettant d'une main une perle en $A$ chaque fois qu'il
  en déposait une autre en $A_2$ avec l'autre main. Après quoi, laissant le
  verre $A$ comme témoin, on verse $A_2$ en un verre $B$ de forme différente.
  Les petits de 4-5 ans concluent alors que la quantité de perles a changé,
  quand bien même ils sont certains que l'on n'a rien enlevé ni ajouté: si le
  verre $B$ est mince et élevé, ils diront qu'il y a <<plus de perles
  qu'avant>> parce que <<c'est plus haut>>, ou qu'il y en a moins parce que
  <<c'est plus mince>>, mais ils s'accorderont à admettre la non-conservation
  du tout.}{piaget:pi}{140}

Piaget postule donc l'existence d'un \emph{stade opératoire} chez
l'enfant, qui est le stade après lequel un enfant maîtrise la notion
d'\emph{opération}.  Dans l'exemple ci-dessus, l'opération est le
transvasage des perles, au cours duquel aucune perle n'est créee ou
détruite; cette opération est réversible, et la notion de conservation
du nombre de perles suit.

Si cette opération n'était pas réversible (par exemple, si les perles
s'aggluttinaient, ou se divisaient, lors du transvasage), elle
pourrait fort bien modifier le nombre de perles présentes dans les
verres, et l'enfant, pour déterminer le nombre de perles dans les
verres, aurait raison de se fier à son intuition.

Piaget pousse plus loin son étude: le nombre $n$ de perles peut être
exprimé comme $1+1+\dots+1$, qu'on voit comme opération <<ajouter $n$
fois une perle dans le verre>>. \fcite{On présente d'emblée au sujet
  les deux verres $A$ et $B$ et on lui demande de mettre simultanément
  une perle dans chaque verre; l'une avec la main gauche, l'autre avec
  la droite; lors de petits nombres (4 ou 5), l'enfant croit d'emblée
  à l'équivalence des deux ensembles, ce qui semble annoncer
  l'opération, mais lorsque les formes changent trop, à mesure que la
  correspondance se poursuit, il renonce à admettre
  l'égalité!}{piaget:pi}{141--142} On voit ainsi que l'enfant accepte
de composer 4 ou 5 fois une opération, mais pas plus.

Le comportement opératoire peut, en première approximation, être
opposé au comportement \emph{intuitif}.  Paradoxalement, une structure
mentale plus rigide --- celle de groupe --- donne plus de souplesse à
un enfant, en particulier lors de l'appréhension de domaines nouveaux.
Piaget parle d'intuition \emph{imagée}, ou \emph{perceptive},
\fcite{tant que les enfants se placent au point de vue de la seule
  perception\dots Mais dès que le sujet renonce à invoquer l'apparence
  sensible pour réfléchir aux transformations comme telles, il est
  conduit à supposer ou à affirmer la conservation. Les opérations qui
  conduisent à ce résultat présentent deux aspects distincts: identité
  et réversibilité.}{piaget:dq}{15--16}; et \fcite{l'opération,
  c'est-à-dire l'action devenue réversible, ne se trouve, à ce niveau
  du développement, qu'en partie dissociée de l'intuition,
  c'est-à-dire de l'action et de la perception
  irréversibles.}{piaget:dq}{17--18} Même l'intuition \emph{articulée}
\fcite{n'est que le produit de régulations successives qui ont fini
  par articuler les rapports globaux inanalysables du début, et non
  pas encore d'un groupement proprement dit.}{piaget:pi}{158}  Pour
Piaget, un \emph{groupement} est une structure abstraite générale,
agissant sur des objets concrets. Il ne s'agit donc pas encore d'un
groupe au sens mathématique; on parlerait plutôt de \emph{magma}, de
\emph{groupoïde}, ou de \emph{relation d'équivalence}. Il n'est pas
impossible que le concept que Piaget recherchait était précisément
celui de catégorie.

Le groupement est l'\fcite{aboutissement d'une suite
  d'équilibrations.}{piaget:pi}{139 bas}  Curieusement, Piaget voit
déjà dans le groupement l'apparition des principes de
conservation (\fcite{\dots tant qu'il n'y a pas groupement, il ne
  saurait y avoir conservation des ensembles ou totalités, tandis
  que l'apparition d'un groupement est attestée par celle d'un
  principe de conservation. Par exemple, le sujet capable de
  raisonnement opératoire à structure de groupement sera
  d'avance assuré qu'un tout se conservera indépendamment de
  l'arrangement de ses parties, tandis qu'il le contestera
  auparavant}{piaget:lela}{245}), alors que la physique classique ne
les voit que dans la structure de groupe (voir \S\ref{ss:physique}).

Un groupe, pour Piaget, est un groupement fermé pour l'opération de
produit et d'inverse. On pense ainsi à un groupement comme à un
ensemble de bijections, i.e.\ d'opérations réversibles; toutefois,
si le groupement est <<trop petit>> au sens où il ne contient pas tous
les produits et inverses de ses éléments, on <<ne sait pas>> que ses
éléments sont inversibles. Ainsi le \fcite{raisonnement
  préconceptuel\dots ne repose que sur des emboîtements incomplets et
  échoue ainsi à toute structure opératoire réversible.}{piaget:pi}{139
  haut} On a toutefois déjà affaire à une catégorie.

Remarquons enfin que Piaget fait très peu mention de l'élément neutre
du groupe dans son discours. Celui-ci pose un problème particulier,
car il est difficile d'accepter une <<transformation>> qui ne
transforme rien, comme le remarque G.  Henriques: \fcite{Il va sans
  dire qu'une transformation identique ne saurait se confondre avec
  l'absence d'activité transformatrice\dots La seule justification
  [des actions transformatrices] réside dans la complétion des
  systèmes de transformations auxquels elles viennent
  s'adjoindre.}{piaget:mc}{194}

Une raison pour la réticence de Piaget peut aussi être que, plus que
la notion de groupe, c'est la notion d'\emph{espace homogène} qui lui
semble essentielle. Un espace homogène est un espace $X$ muni d'une
\emph{action} d'un groupe $G$, i.e.\ d'une opération attribuant à
$g\in G$ et $x\in X$ un élément $gx\in X$, celle-ci satisfaisant les
axiomes $1x=x$ et $g(hx)=(gh)x$. Un groupe est un cas particulier
d'espace homogène; on prend $X=G$ et on définit l'action comme la
multiplication du groupe. Une différence essentielle entre un groupe et
un espace homogène est que ce dernier n'a, au contraire de $G$, aucun
point privilégié.  C'est donc un pas de plus dans la direction de la
<<décentration>>, de la capacité de percevoir comme un acteur de
l'univers monde au même titre que d'autres, et non comme un point
central de celui-ci. La \fcite{décentration ou coordination des
  rapports successivement centrés constitue par conséquent un
  <<groupement>> opératoire}{piaget:mv}{171}; et \fcite{le groupement
  est le contraire de l'égocentrisme.}{piaget:dq}{275}

Piaget voit dans la décentration un des pas essentiels de l'abstraction:
\fcite{Mais cet état initial, qui se retrouve en chacun des domaines de la
  pensée intuitive, est progressivement corrigé grâce à un système de
  régulations, qui annoncent les opérations. Dominée d'abord par le rapport
  immédiat entre le phénomène et le point de vue du sujet, l'intuition évolue
  dans le sens de la décentration. Chaque déformation poussée à l'extrême
  entraîne la réintervention des rapports négligés. Chaque mise en relation
  favorise la possibilité d'un retour. Chaque détour aboutit à des
  interférences qui enrichissent les points de vue. Toute décentration d'une
  intuition se traduit en une régulation, qui tend dans la direction de la
  réversibilité, de la composition transitive et de l'associativité, donc, au
  total, de la conservation par coordination des points de
  vue.}{piaget:pi}{148--149}

\subsection{Groupes de déplacements}\label{ss:gd}
On peut voir l'espace tridimensionel, dans lequel nous vivons, comme
le groupe de ses translations. Ce faisant, on fait l'amalgame entre le
groupe (les \emph{translations} de l'espace $\R^3$) et l'espace
homogène (de nouveau $\R^3$). On peut toutefois revenir au point de
départ dans l'espace, mais dans une autre orientation; ceci se
décrit naturellement par un groupe de déplacements plus riche
que celui des translations. Du point de vue pratique, la maîtrise
de ce groupe est essentielle pour l'enfant s'il souhaite se
déplacer efficacement.  Comme le note Piaget, il y a \emph{groupe}
lorsqu'il y a \fcite{un ensemble coordonné de déplacements
  susceptibles de revenir à leur point de départ, et tel que
  l'état final ne dépende pas du chemin
  parcouru.}{piaget:cr}{95} Si l'orientation de l'objet peut changer
au cours du déplacement, il faut que le groupe en tienne
compte.\footnote{Il existe un sous-groupe du groupe des
  déplacements, constitué des différentes orientations que
  peut prendre un objet après son retour à l'origine. Il s'agit
  du \emph{groupe d'holonomie}, dont nous n'aborderons pas l'étude
  psychologique, bien qu'elle soit certainement très prometteuse.}

Piaget distingue six stades de développement de groupes de
déplacement chez l'enfant. Au début, les groupes ne sont
présents que comme <<réceptacles>> du système de
coordonnées de l'univers, et \fcite{on ne saurait parler de
  groupes, impliquant la perception de rapports actifs que le sujet
  établit entre les choses et lui.}{piaget:cr}{100} Ils sont donc d'une
part entièrement concrets, et d'autre part aucune opération ne
leur est appliquée. Ceci caractérise les deux premiers stades.

Les groupements apparaissent au troisième stade. Ils sont toujours concrets, et
ne servent qu'à coordonner des domaines qui jusque-là étaient disjoints:
\fcite{Les opérations dont il s'agit ici sont donc des <<opérations concrètes>>
  et non pas encore formelles: toujours liées à l'action, elles structurent
  celle-ci logiquement\dots mais n'impliquent en rien la possibilité de
  construire un discours logique indépendant de l'action.}{piaget:pi}{155} On
ne voit pour l'instant qu'un début d'objectivisation de l'espace, i.e.\ la
\fcite{constitution d'un champ indépendant du corps propre et dans lequel
  celui-ci se déplace comme un objet parmi d'autres objets.}{piaget:cr}{107}

Au quatrième stade, il y a dissociation des déplacements de l'objet et de ceux
du sujet: l'enfant a conscience à la fois de ses propres déplacements et de
l'immobilité des objets, même si \fcite{le sujet demeure égocentrique,
  géométriquement parlant: il ne conçoit pas encore les positions et
  déplacements comme relatifs les uns aux autres, mais uniquement comme
  relatifs à lui.}{piaget:cr}{160} La réversibilité apparaît aussi à ce stade.

Le cinquième stade voit l'apparition de groupes proprement objectifs. C'est ici
que l'enfant acquiert la notion de \fcite{déplacement des objets les uns par
  rapport aux autres}{piaget:cr}{160 bas}; autrement dit, il parvient à
percevoir la relation entre deux objets (et non seulement sa relation à un
objet) à l'aide d'un groupe. Il peut aussi tenir compte de déplacements
successifs d'une chose qu'il cherche: \fcite{auparavant, il cherchait
  systématiquement l'objet là où il l'avait trouvé une première
  fois.}{piaget:cr}{161}

C'est enfin au sixième stade que l'enfant perçoit et utilise véritablement le
groupe. Il peut en particulier utiliser des <<détours>>, i.e.\ par exemple
contourner un objet difficile à franchir comme un lit. En termes de groupes, ce
sont des équations du genre $w=xyz$ qu'il faut utiliser, où $w$ est <<traverser
le lit par le milieu>> et $x,y,z$ sont respectivement <<aller au coin du lit>>,
<<suivre son bord>> et <<retourner au milieu de l'autre côté>>. \fcite{L'espace
  est constitué pour la première fois à titre de milieu immobile. Cette
  acquisition finale garantit ainsi l'objectivité des groupes perçus et la
  possibilité d'étendre ces groupes aux déplacements ne tombant pas directement
  dans le champ de la perception.}{piaget:cr}{181}

Les groupes apparaissent dans le cadre des déplacements. Ils sont ensuite
présents dans d'autres contextes, en particulier la pensée formelle.

\subsection{Réversibilité}\label{ss:rev}
C'est dans la réversibilité que Piaget voit le développement de l'abstraction:
\fcite{en atteignant le niveau de la réversibilité entière, les opérations
  concrètes issues des régulations précédentes se cordonnent en effet en
  structures définies\dots qui se conservent la vie durant, sans exclure la
  formation de systèmes supérieurs.}{piaget:lela}{217}

La pensée formelle nécessite la possibilité de spéculations \emph{réversibles}
sur l'univers: \fcite{avec la pensée formelle\dots une inversion de sens
  s'opère entre le réel et le possible. Au lieu que le possible se manifeste
  simplement sous la forme d'un prolongement du réel ou des actions exécutées
  sur la réalité, c'est au contraire le réel qui se subordonne au possible: les
  faits sont dorénavant conçus comme le secteur des réalisations effectives au
  sein d'un univers de transformations possibles.}{piaget:lela}{220} Il est
important, lors de telles spéculations, de pouvoir facilement modifier de façon
réversible un état hypothétique du monde, comme un joueur d'échecs déplace
aisément mentalement une pièce dans ce qu'il considère comme pouvoir être la
position de l'échiquier dix coups plus tard.

La réversibilité est utilisée par Piaget comme critère de constitution des
opérations concrètes, et c'est le recours à celle-ci qui permet la constitution
de notions de conservation. La réversibilité est l'argument utilisé par les
enfants pour justifier, par exemple, la conservation de la pâte à modeler que
l'on transforme en boulettes: \fcite{C'est précisément la réversibilité qui
  engendre la conservation\dots la réversibilité est le processus même dont la
  conservation est le résultat\dots l'identité, en fait, devient un argument à
  partir du moment où elle est subordonnée à la réversibilité\dots C'est donc
  bien la réversibilité qui entraîne la conservation et non pas
  l'inverse.}{piaget:bc}{293}

La réversibilité garantit aussi la possibilité de l'équilibre: de même qu'il y
a équilibre de pression entre deux chambres si et seulement s'il y a
possibilité de flot, dans les deux directions, d'air entre les chambres,
\fcite{tout état d'équilibre est reconnaissable\dots à une certaine forme de
  réversibilité.}{piaget:lela}{240}

\section{Les groupes en et hors des mathématiques}\label{s:math}
La théorie des groupes est, bien évidemment, apparue en mathématiques. Il est
important, toutefois, de remarquer qu'elle n'a pas été créée pour elle-même,
mais pour expliquer d'autres parties de mathématiques: la résolution
d'équations polynomiales (Galois) ou différentielles (Lie). Son étude
systématique n'a commencé qu'un demi-siècle plus tard.

Piaget a reconnu très tôt l'application possible de la théorie des
groupes à la psychologie. Reprenons l'exemple du comptage de perles
cité au chapitre~\ref{s:gpespiaget}: un équilibre est atteint quand l'enfant
comprend que les variations des hauteur et de largeur des perles dans
le verre $B$ non seulement peuvent, mais \emph{doivent} se compenser.
Piaget décrit cet équilibre mobile par les caractéristiques
suivantes:\fcite{\begin{enumerate}
  \item deux actions successives peuvent se coordonner en une seule;
  \item le schème d'action, déjà à l'oeuvre dans la pensée intuitive, devient
    réversible;
  \item un même point peut être atteint, sans être altéré, par deux voies
    différentes;
  \item le retour au point de départ permet de retrouver celui-ci identique à
    lui-même;
  \item la même action, en se répétant, ou bien n'ajoute rien à elle-même, ou
    bien est une nouvelle action avec effet cumulatif.
  \end{enumerate}
  On reconnaît bien là la composition transitive, la réversibilité,
  l'associativité et l'identité}{piaget:pi}{155}

Au contraire, la commutativité ($ab=ba$) ne semble pas être une notion aussi
importante, soit pour décrire des phénomènes physiques, soit chez l'enfant. Il
suit ainsi que la règle d'inversion d'un produit, $(ab)^{-1}=b^{-1}a^{-1}$, est
comprise par l'enfant dès qu'il maîtrise la notion d'inverse, comme le montre
l'expérience suivante: \fcite{On enfile le long d'un même fil de fer trois
  bonshommes de couleurs différentes $A,B,C$, ou bien l'on fait entrer dans un
  tube de carton (sans chevauchements possibles) trois boules $A,B,C$. On fait
  dessiner à l'enfant le tout, à titre d'aide-mémoire. Puis on fait passer les
  éléments $A,B,C$ derrière un écran ou à travers le tube et l'on fait prévoir
  l'ordre direct de sortie (à l'autre extrémité) et l'ordre inverse de retour.
  L'ordre direct est prévu par tous. L'ordre inverse, par contre, n'est acquis
  que vers 4-5 ans, à la fin de la période préconceptuelle.}{piaget:pi}{145}

C'est cette même approche --- reconnaître une structure de groupe sous-jacente
--- qui a été appliquée innombrablement dans divers domaines scientifiques.
Nous décrivons ci-dessous quelques telles situations où des groupes
apparaissent naturellement pour décrire la réversibilité d'opérations dans un
univers.

\subsection{Mathématiques}
Les premiers exemples de groupes rencontrés dans un cours de
mathématiques sont généralement les groupes numériques,
donnés par l'addition d'entiers, ou l'addition <<modulo $n$>>, ou
encore la multiplication de rationnels non-nuls.  C'est ce que Piaget
décrit comme les <<groupements additifs>> et les <<sériations
d'ordre>> (si l'on considère la relation <<plus petit que>> sur les
rationnels).

L'arithmétique sur les entiers se propose d'étudier les
propriétés du groupe $\mathbb Z$ des entiers relatifs. Les
mathématiques construisent $\mathbb Z$ à partir des entiers
naturels $\mathbb N$, et construisent ces derniers à partir des
ensembles, réduisant ainsi l'arithmétique à la logique.
Piaget constate que l'enfant suit le même chemin dans son
développement: \fcite{les opérations constitutives du
  nombre\dots ne requiert pas autre chose que les groupements additifs
  de l'emboîtement des classes\dots mais fondu en un seul tout
  opératoire\dots\ Cette réunion de la différence et de
  l'équivalence suppose\dots le passage du logique au
  mathématique.}{piaget:pi}{154}

Toutefois, une source inépuisable, et plus générale, de groupes est la
suivante. On prend \emph{n'importe quel} objet $X$, muni de ses
propriétés; par exemple, une figure et l'alignement de ses parties; un
ensemble de points et leurs distances respectives, un graphe et
l'ensemble de ses arêtes. On considère l'ensemble $G$ de toutes les
applications bijectives\footnote{c'est-à-dire mettant les objets au
  source et au but en correspondance bi-univoque.} $f:X\to X$ qui
préservent ces propriétés.  C'est ce qu'on appelle le \emph{groupe
  d'automorphismes} de $X$. Ce groupe est un objet algébrique dont,
d'une part, l'étude peut révéler de l'information sur $X$, permettre
de distinguer $X$ d'un $X'$ \emph{à priori} similaire, etc; et d'autre
part, dont la simple existence révèle des symétries de $X$, permettant
de simplifier grandement les calculs effectués sur $X$.

\subsection{Le groupe INRC}
Ce groupe apparaît à de nombreuses reprises dans l'oeuvre de Piaget; c'est un
groupe de symétries d'énoncés logiques. Il s'agit donc d'une abstraction au
deuxième degré, le premier degré étant la constitution d'un énoncé logique à
partir d'une situation concrète.

Bien qu'il ne semble pas aussi important que les autres groupes apparaissant
dans ce travail, décrivons-le brièvement. Les lettres I,N,R,C indiquent
respectivement les transformations \emph{identité}, \emph{négation},
\emph{réciproque} et \emph{corrélative}, i.e.\ à partir de l'énoncé <<$A$
implique $B$>>, les énoncés <<$A$ implique $B$>>, <<$B$ implique $A$>>,<<$A$
n'implique pas $B$>>, et <<$B$ n'implique pas $A$>>.

On a donc affaire à un groupe de transformations d'énoncés logiques, à quatre
éléments. Il est commutatif, tous ses éléments sont leurs propres inverses, et
il se nomme généralement le \emph{groupe de Klein}.

Si on généralise la construction d'INRC aux symétries
d'énoncés logiques avec plus que deux variables, on obtient
toujours un groupe abélien, mais plus gros. On a par exemple un
groupe d'ordre $8$ pour décrire les symétries de formules
logiques à deux arguments et un résultat; et d'autres groupes,
toujours isomorphes à INRC, dans des systèmes logiques plus
complexes, comme par exemple la logique où toutes les opérations
sont ternaires~\cite{piaget:tol}*{page 159}.

L'intérêt d'un groupe, dans un tel contexte, est de préserver dans une
structure relativement simple des changements opérés sur une formule. Plutôt
que de réécrire, ou de reconsidérer mentalement, une formule à laquelle on a
appliqué un élément d'INRC, on peut conserver la formule telle quelle et
considérer l'élément d'INRC en question. C'est aussi cette
philosophie d'économie que l'on retrouvera plus bas, au \S\ref{ss:info}.

\subsection{Physique}\label{ss:physique}
La physique se préoccupe de la construction d'\emph{invariants},
i.e.\ de grandeurs, ou propriétés d'un système qui ne
dépendent pas de choix arbitraires faits pour effectuer des
mesures. Ainsi la position d'un objet dépend d'un choix de
référentiel; mais son énergie n'en dépend pas. De même
toute grandeur physique mesurable est \emph{par nécessité} un
invariant --- invariant du groupe des changements de
référentiel.

S'il existe plusieurs choix de référentiels, tous équivalents entre eux, alors
il existe un groupe de transformations convertissant ces référentiels les uns
en les autres.  Il ressort \emph{in fine} que ce qui est invariant est
précisément ce qui est \emph{conservé} par un groupe, et c'est là l'une des
idées essentielles de la mécanique classique. Généralisant, \fcite{\dots la
  conservation constitue une condition nécessaire de toute activité
  rationnelle.}{piaget:gn}{6}

Piaget était intimement conscient de ce rapport entre invariants et groupes:
\fcite{La découverte d'une notion de conservation par l'enfant est toujours
  l'expression de la construction d'un\dots <<groupe>> (mathématique)
  d'opérations.}{piaget:gs}{354}, \fcite{Les quantités qui se conservent sont
  des invariants de groupements ou de groupes}{piaget:dq}{214}

On peut se demander lequel, de l'invariant ou du groupement, est construit en
premier par l'enfant. Piaget écrit: \fcite{La condition commune à tous les
  groupements [est de] définir les transformations en fonction d'invariants et
  réciproquement.}{piaget:dq}{74} Toutefois, on lit aussi que \fcite{là où il y
  a <<groupement>> il y a conservation d'un tout, et cette conservation
  elle-même ne sera pas simplement supposée par le sujet à titre d'induction
  probable, mais affirmée par lui comme une certitude de sa
  pensée.}{piaget:pi}{150} On voit donc que la construction d'un groupe confère
de la solidité à un invariant, et peut être considérée comme préalable.

Comme on l'a vu au \S\ref{ss:rev}, c'est via la réversibilité que l'enfant
accède à la notion de groupe puis d'invariant: dans l'expérience des perles
décrite au chapitre~\ref{s:gpespiaget}, le progrès de l'enfant apparaît quand
\fcite{il répond qu'un transvasement de $A$ en $B$ peut être corrigé par le
  transvasement inverse}{piaget:pi}{150}

Il est intéressant d'examiner de plus près l'utilisation réelle des groupes en
physique. Un résultat essentiel est communément appelé le \emph{théorème de
  Noether}. Brièvement, il dit qu'à tout groupe de transformations géométriques
(translations, rotations, etc.) laissant invariantes les lois physiques
correspond une loi de conservation. Ainsi la conservation de l'énergie est liée
à l'invariance par décalage temporel; la conservation de l'impulsion est liée à
l'invariance par translations; et la conservation du moment angulaire est liée
à l'invariance par rotations des lois de la physique.

On peut exprimer plus précisément ces résultats. Considérons d'abord le
formalisme galiléen; on a alors un système physique, par exemple une particule,
qui évolue dans un univers. Le système est décrit comme un point dans espace,
dans notre exemple $M=\R^3\times\R$, où $\R^3$ désigne la position dans
l'espace et $\R$ est la coordonnée <<temps>>.  Il y a des choix arbitraires: la
position et l'orientation d'un repère dans $\R^3$, et l'origine du temps. On
considère donc un groupe, le \emph{groupe galiléen}, qui agit sur l'espace $M$.
Il est engendré par les transformations de la forme suivante:
\begin{itemize}
\item $(x,t)\mapsto(x+vt,t)$, le mouvement rectiligne uniforme, pour tout $v\in\R^3$,
\item $(x,t)\mapsto(x+p,t+s)$, les translations, pour tous $p\in\R^3$ et $s\in\R$, et
\item $(x,t)\mapsto(Gx,t)$, les rotations, pour tout $G\in O(3)$ (matrice de
  rotation).
\end{itemize}
C'est un groupe de dimension $10$, et la mécanique classique
postule que les lois de la physique doivent être invariantes par ce
groupe; ainsi il n'y a pas de point ni de direction privilégiée,
ni dans le temps ni dans l'espace, et les lois ne dépendent que des
vitesses relatives.

La mécanique relativiste, aussi dite de Lorentz, ne donne pas de rôle
singulier à la coordonnée <<temps>>; elle considère $M$ comme $\R^4$,
mais avec une nouvelle notion de distance, qui peut parfois prendre
des valeurs négatives. Le groupe d'invariants est alors le
\emph{groupe de Lorentz}, engendré par les translations: $\R^4$, et
les rotations de l'espace-temps: $O(3,1)$. C'est à nouveau un groupe
de dimension $10$, mais il semble, à posteriori, mathématiquement plus
naturel que le groupe de Galilée.

Venons-en au formalisme lagrangien; pour référence on
suit~\cite{arnold:cm}*{pages 53 et suivantes}. On a un espace de
configurations $M$, qui est une variété différentielle, par
exemple $\R^3$ pour l'espace des configurations d'une particule dans
l'espace. On ne représente plus une particule comme un chemin
$x(t)$ dans $M$, mais plutôt comme un chemin $(x(t),\dot x(t))$ dans
l'espace tangent $TM$ de $M$; ces points de vue sont évidemment
équivalents.

Un \emph{système mécanique lagrangien} est donné par cet
espace $M$ de configurations, et par une fonction $L$, dite
lagrangienne, sur $TM$. Alors on décrète que l'évolution du
système ne peut se produire que selon des courbes
\emph{extrémales} pour $L$: si $x(t)$, pour $t_0\le t\le t_1$, est
une évolution du système, alors la fonctionnelle
\[\Phi(x)=\int_{t_0}^{t_1}L(x(t),\dot x(t))\]
doit prendre une valeur extrémale en $x$.  En dérivant sous
l'intégrale, puis en intégrant par parties, on trouve que
\[\Phi(x)\text{ est extrémale si et seulement si }\frac{\partial
  L}{\partial x}-\frac{d}{dt}\frac{\partial L}{\partial\dot x}=0.\]
La fonction langrangienne indique ainsi la <<direction de changement
maximal>> sur l'espace tangent.

Dans ce langage, $M$ a naturellement un groupe de difféomorphismes,
qui remplace le groupe de Galilée ou de Lorentz. Le théorème
de Noether formalise de la façon précise suivante l'idée que
toute loi de conservation est justiciable d'un invariant de groupe:
<<si la fonction lagrangienne est invariante par un groupe à un
paramètre $h^s$, pour $s\in\R$, i.e. si $\{h_s\}$ est un groupe de
difféomorphismes de $M$ isomorphe à $\R$, et tel que
$\frac{d}{ds}L(h^sx,h^s\dot x)=0$, alors la quantité $C(x,\dot
x)=\frac{\partial L}{\partial\dot x}\frac{dh^s(x)}{ds}|_{s=0}$ est
conservée.>>

Si on prend comme ci-dessus $M$ l'espace euclidien, et on définit
$L$ comme la différence entre l'énergie cinétique et
l'énergie potentielle, on retrouve les lois de conservation de
Galilée. Par exemple, imaginons pour simplifier que la particule,
qui a une masse $m$, ne peut se déplacer que verticalement, et est
une hauteur $h$. Elle a donc une énergie potentielle $E_p=mgh$
où
$g$ est la constante de gravitation ($9.81ms^{-2}$ au sol) et une
énergie cinétique $E_c=\frac m2(\dot h)^2$, d'où
$L=E_p-E_c=mgh-\frac m2(\dot h)^2$. Si $\Phi(h)$ est extrémale, on
a
\[\frac{\partial L}{\partial h}-\frac{d}{dt}\frac{\partial
  L}{\partial\dot h}=0\text{ si et seulement si }mg=m\ddot h,\]
c'est-à-dire que l'accélération de la particule est
constante, égale à $g$.

Notons finalement que la mécanique quantique fait un usage encore bien plus
étendu de la théorie des groupes. D'une part, l'évolution d'un système dans le
temps est décrit par une transformation inversible (unitaire, plus
précisément); le groupe des difféomorphismes de l'espace des configurations est
de dimension infinie. Les particules élémentaires apparaissent alors
naturellement comme représentations de ce groupe; à chaque représentation
irréductible correspond un type de particule élémentaire. Nous ne nous
aventurerons pas plus dans cette direction.

\subsection{Informatique}\label{ss:info}
Contrairement à une idée communément répandue, l'informatique est liée
à la théorie des groupes par plusieurs liens très forts. Le premier
est certainement la notion de <<pile>>: une structure dans laquelle on
peut stocker des objets, et les en sortir dans l'ordre inverse
d'entrée. La pile est un objet, muni d'une transformation inversible
<<empiler>>: si elle n'était pas inversible, il ne serait pas possible
de récupérer l'objet empilé précédemment. L'expérience de Piaget,
décrite au chapitre~\ref{s:math}, illustre bien le concept de pile.
Elle continue ainsi, par des manipulations de cette pile:
\textsf{Après quoi on imprime un mouvement de rotation de $180^o$ à
  l'ensemble du dispositif (fil de fer ou tube) et l'on fait prévoir
  l'ordre de sortie (qui est donc renversé). L'enfant ayant contrôlé
  lui-même le résultat, on recommence, puis on effectue deux
  demi-rotations ($360^o$ en tout), puis trois, etc.}
  
\fcite{Or, cette épreuve permet de suivre pas à pas tous les progrès de
  l'intuition jusqu'à la naissance de l'opération. De 4 à 7 ans, le sujet
  commence par ne pas prévoir qu'une demi-rotation changera l'ordre $ABC$ en
  $CBA$; puis, l'ayant constaté, il admet que deux demi-rotations donneront
  aussi $CBA$.  Détrompé par l'expérience, il ne sait plus prévoir l'effet de
  trois demi-rotations.}{piaget:pi}{145}

On peut imaginer, à titre d'exemple, la programmation de la
commande <<undo>> (défaire) dans un programme d'édition de
texte. Au premier jet, un programmeur construira une pile, dans
laquelle il insérera, au fur et à mesure de la modification du
document, un historique de l'état du programme après chaque
commande. Dans un deuxième jet, il créera une fonction <<undo>>
qui dépend du type de la commande exécutée. Le problème de
cette approche est que si le programme est complexe, les
différentes combinaisons de <<undo>> le sont encore plus; de plus,
chaque commande nécessite un espace de stockage complètement
différent des données détruites.

La solution habituellement adoptée est celle du <<pattern
commande>>, décrite par exemple dans l'ouvrage classique de
programmation dite orientée-objet~\cite{gamma-:dp}. L'idée est
d'\fcite{encapsuler une requête comme un objet, autorisant ainsi le
  paramétrage des clients par différentes requêtes, files
  d'attente et récapitulatifs de requêtes, et de plus,
  permettant la réversion des opérations.}{gamma-:dp}{271} Dans
notre exemple de traitement de texte, toutes les commandes seraient
implémentées comme des \emph{objets}. En tant que tels, ils
peuvent être assemblés en \emph{séquences}; ceci revient à
les multiplier dans une structure algébrique. Ils sont tous munis
de leur propre opération <<undo>>, qui défait l'effet
précédent de la commande. C'est précisément l'inverse
qu'on avait dans une structure de groupe.

Ainsi l'introduction de la notion d'opération réversible
revêt une importance capitale du point de vue de l'économie de
moyens. Même avec le développement galopant de la technologie,
il est toujours d'actualité d'utiliser le minimum de mémoire
nécessaire dans un programme; et si ce critère est valable pour
les ordinateurs, il l'est à fortiori pour les êtres humains.

Si on admet qu'on manipule un objet --- en l'occurrence, un fichier de
texte --- de manière irréversible, alors pour garantir la
possibilité de revenir en arrière on doit stocker toutes les
étapes intermédiaires de cet objet. Si au contraire on admet que les
transformations sont réversibles --- par exemple, si l'objet texte
contient en lui-même les moyens de défaire les modifications
qu'il a subies --- alors on ne doit stocker qu'une seule copie de
l'objet.

L'exemple du joueur d'échecs ne se souvenant que d'une position de
l'échiquier à la fois, et manipulant les pièces sur cet échiquier de
façon réversible, a déjà été évoqué au \S\ref{ss:rev}. Les
programmes informatiques d'échecs fonctionnent d'ailleurs selon ce
même principe.  Finalement, il est intéressant de noter que les études
du psychologue des échecs Krogious montrent qu'une des sources
majeures d'erreur de calcul de joueurs d'échecs vient précisément de
mouvements mentaux modifiant irréversiblement la nature de
l'échiquier, c'est-à-dire essentiellement lors de captures de
pièces~\cite{krogious:psy}.

Remarquons que, du point de vue théorique, on peut toujours
supposer que les calculs effectués par un ordinateur le sont de
manière réversible. C'est particulièrement intéressant
dans le contexte des ordinateurs quantiques, qui sont soumis, comme
toute la physique quantique, à n'exécuter que des opérations
réversibles. Même pour un ordinateur classique, ce résultat
présente un intérêt pratique: toute opération
non-réversible entraîne un accroissement d'entropie, donc une
dissipation d'énergie sous forme de chaleur. Les processeurs
actuels sont limités en particulier par ce problème
d'évacuation de la chaleur; or un ordinateur classique, mais
opérant de manière réversible, n'est pas soumis à cette
limite~\cite{bennett:reversibility}.

Ce résultat de réversibilité peut prendre la forme précise
suivante: soit $M$ une machine de Turing\footnote{C'est un modèle
  canonique d'ordinateur élémentaire}. Alors il existe une
machine de Turing équivalente $M'$ telle que si $U\to V$ est une
instruction de $M'$, alors $V\to U$ est aussi une instruction de $M'$.
Il y a de plus des bornes sur la complexité de $M'$ en fonction de
la complexité de $M$~\cite{sapir-:isom}*{Lemme~3.1}.

\section{Jeux éducatifs}
Les travaux de Piaget ont eu une grande influence sur le domaine de
l'éducation, la plus importante étant peut-être d'avoir mis en
évidence que les grands pas effectués par un enfant au cours de son
développement intellectuel sont autant favorisés, voire plus, par ses
pratiques ludiques que par les explications rationnelles de ses
professeurs. En cela, il s'est opposé à la pédagogie traditionelle qui
\fcite{considérait le jeu comme une sorte de déchet mental ou tout au
  moins comme une pseudo-activité, sans signification fonctionnelle et
  même nuisible aux enfants, qu'il détourne de leurs devoirs\dots\
  [Puis, citant K. Gross:] Le jeu est <<préexercice>>, et non pas
  seulement exercice, parce qu'il contribue au développement de
  fonctions dont l'état de maturité n'est atteint qu'à la fin de
  l'enfance.}{piaget:fs}{158} Ce dernier chapitre a pour but de
décrire des jeux éducatifs naturels dans lesquels les concepts
fondamentaux de la théorie des groupes peuvent se trouver. On
indiquera aussi différentes pistes pour l'implémentation informatique
de ces jeux.

Insistons encore sur le fait qu'aucun concept de théorie des groupes
ne doit être apparent dans les jeux ci-dessous; toute cette structure
doit être sous-jacente, si l'enfant est appelé à \emph{découvrir} et
non \emph{apprendre}. Ainsi, de la même manière que Piaget voulait
\fcite{trouver un biais de manière à éviter toute allusion aux
  désignations linguistiques de la multiplication}{piaget:ar:cm}{31}
dans ses expériences, dans le but d'éviter que les sujets répondent ce
qu'ils ont appris <<par coeur>> sans le comprendre, on évitera de
désigner explicitement le groupe utilisé: Piaget préfèrera, par
exemple, ``J'ai 21 crayons et 15 gommes. Je veux donner le même nombre
de crayons et de gommes à chacun. Combien d'amis puis-je inviter au
maximum'' à ``Quel est le p.g.c.d. de 21 et 15?''.  Il s'agit d'amener
l'enfant à \fcite{reconnaître en en un objet un caractère $x$ pour
  l'utiliser à titre d'élement d'une structure différente de celle des
  perceptions considérées, ce que nous désignerons alors du nom
  d'abstractions et de généralisation
  <<constructive>>}{piaget:mp}{395}.

\subsection{Principe général}\label{ss:factor}
La donnée d'un groupe permet de construire assez facilement un jeu
classique qu'on peut appeler \emph{jeu de factorisation}.

Dans l'enseignement mathématique de niveau secondaire, par
factorisation on entend généralement la factorisation ou décomposition
d'un nombre entier en produit de nombres premiers --- un nombre
premier étant un nombre qui est indécomposable en produit de nombres
plus petits.  Par exemple, si le nombre à factoriser est 84, la
solution du \emph{jeu de factorisation} est $84=2\cdot2\cdot3\cdot7$.

La généralisation de ce jeu à un groupe quelconque est la suivante: on
fournit au joueur un élément du groupe --- dans notre exemple, 84 ---
et un \emph{ensemble générateur} du groupe, c'est-à-dire des
éléments du groupe qui ont la propriété que n'importe quel autre
élément du groupe peut s'écrire comme un produit d'un nombre fini de
ces générateurs. Dans notre exemple, ces générateurs sont tous les
nombres premiers, et la propriété qui veut que tout nombre entier
puisse s'écrire comme un produit de nombres premiers est un théorème
fondamental de la théorie des nombres.\footnote{L'exemple donné
  concerne l'ensemble des nombres entiers positifs avec la
  multiplication, qui forment en réalité un monoïde (voir le
  \S\ref{ss:cat}) et non un groupe. On pourrait toutefois
  formuler cet exemple de manière tout à fait équivalente dans le
  groupe des nombres rationnels --- les fractions de nombres entiers
  --- non nuls, où cette fois tous les éléments ont un inverse pour la
  multiplication.}

Il est possible de débarrasser ce jeu de son apparence très
mathématique, en choisissant à la place des nombres entiers un
\emph{groupe concret}, tel qu'il a été présenté au
chapitre~\ref{s:intromath}. Les éléments d'un groupe concret $G$ sont
des bijections d'un ensemble $U$, qui sera généralement judicieusement
choisi pour avoir une représentation visuelle. Il est important de
réaliser que seuls les éléments de $U$ sont des entités visibles ou
palpables. Le groupe $G$ ne contient que des opérations
(transformations) sur les éléments de $U$.

Pour la description des jeux qui suivront, il est utile d'introduire
une nouvelle notation des groupes concrets, qui est celle dite de
l'\emph{action de groupe}\footnote{Pour les mathématiciens: on ne
  parlera ici que d'action fidèle.}: si $G$ est un groupe concret sur
un ensemble $U$, on dit que $G$ \emph{agit sur $U$}, et pour une
bijection $g\in G$ et un élément $u\in U$, on note par $g\cdot
u=g(u)\in U$ le résultat de la bijection $g$ appliquée à l'élément
$u$.  Cette notation est conforme à l'intuition grâce aux propriétés
suivantes:
\begin{enumerate}
\item $(g_1g_2)\cdot u=g_1\cdot(g_2\cdot u)$, pour tous $g_1,g_2\in G$
  et pour tout $u\in U$;
\item $1\cdot u=u$ pour tout $u\in U$, où $1$ désigne l'élément neutre
  du groupe, qui n'est autre que la transformation identique: celle
  qui opère sur chaque élément de $U$ en le laissant inchangé.
\end{enumerate}

Dans sa forme classique, le jeu de factorisation demande une propriété
supplémentaire à l'action de groupe, qui est appelée
\emph{transitivité}: l'action de $G$ sur $U$ est dite
\emph{transitive} s'il est toujours possible de passer d'un élément de
$U$ à un autre en utilisant un opération du groupe; autrement dit, si
pour tous $u_1,u_2\in G$ il existe une opération $g\in G$ telle que
$g\cdot u_1=u_2$.  Lorsque cette propriété de transitivité n'est pas
vérifiée, cela signifie qu'un point donné $u\in U$ ne peut être
transformé qu'en certains autres points de $U$. L'ensemble des points
atteignables à partir de $u$ est appelé l'\emph{orbite} de
$u$\footnote{Par analogie bien sûr à l'orbite d'une planète, qui
  représente tous les points de l'espace pouvant être atteints par la
  planète sous l'<<action>> du temps.}, et l'ensemble $U$ se retrouve
ainsi stratifié en orbites disjointes.

En remplaçant alors $U$ par une de ces orbites $O$, on se retrouve
avec le même groupe $G$, qui agit cette fois transitivement sur $O$.
On pourrait donc penser que la situation est fondamentalement la même
que celle d'une action transitive; on verra toutefois plus loin qu'il
est possible de construire des jeux très intéressants sur des groupes
concrets avec une action non transitive.

Cette nouvelle notation permet d'exprimer très naturellement les
données de notre jeu de factorisation: si $G$ est un groupe concret
agissant transitivement sur un ensemble $U$, si $u_0,u_f$ sont deux
éléments de $U$ et si $S$ est un ensemble générateur de $G$, le
but du jeu est de trouver des éléments (pas forcément distincts)
$s_1,s_2,\dots,s_n\in S$ tels que $s_n\cdot s_{n-1}\cdots s_1\cdot
u_0=u_f$; le fait que l'action soit transitive assure que ce jeu
possède toujours une solution.

Remarquons qu'on a écrit le produit $s_n\cdot s_{n-1}\cdots s_1\cdot
u_0$, alors que $s_1$ agit en premier, $s_2$ en second, et $s_n$ en
dernier. C'est dû à notre convention d'agir par la gauche sur les
éléments. Si l'on agissait par la droite, on trouverait les éléments
$s_1,\dots, s_n$ dans leur ordre naturel.

\subsection{Exemples classiques}\label{ss:class}
\subsubsection{Orientation sur un plan quadrillé}
Dans beaucoup des jeux des premiers temps de l'informatique, le joueur
se retrouvait sur un écran quadrillé, avec la possibilité de s'y
déplacer à l'aide des quatre touches de direction, que nous
appellerons Sud, Ouest, Nord et Est. Ces quatre touches représentent
un ensemble générateur du groupe des déplacements sur un
quadrillage. Plus précisément, l'ensemble $U$ est l'ensemble des
points du quadrillage et le groupe $G$ est l'ensemble des couples de
nombres entiers permettant de représenter les déplacements, par
exemple avec la convention que $(3,-2)$ signifie <<trois cases à
l'est, moins deux cases au nord, c'est à dire deux cases au sud>>.
Dans cette convention, on écrirait Sud=$(0,-1)$, Ouest=$(-1,0)$,
Nord=$(0,1)$, E=$(1,0)$, l'ensemble \{Sud,Ouest,Nord,Est\} étant un
ensemble générateur pour le groupe $G$.

Le jeu de factorisation consiste ici simplement de se rendre d'un
point à un autre en utilisant les quatre flèches. Cet exemple montre
donc qu'un jeu basé sur la théorie des groupes peut être très simple
lorsque le groupe sous-jacent l'est aussi.

\subsubsection{Le cube hongrois}
Ce casse-tête bien connu, dû au professeur Ernö Rubik, est un cube
qu'on voit comme composé de petits cubes à raison de trois pour la
hauteur, trois pour la largeur et trois pour la profondeur, soit
vingt-six au total, le cube central n'étant pas visible. Au départ,
tous les petits carrés visibles sur une même face du grand cube sont
coloriés de la même couleur. Le mécanisme de l'objet permet de faire
tourner n'importe quelle face du grand cube sur elle-même,
mélangeant ainsi les couleurs des faces contiguës.

Le jeu de factorisation consiste à partir d'une configuration
quelconque du cube, et de retrouver la configuration de départ en
utilisant les mouvement de rotations des faces (l'ensemble
générateur). Par opposition au jeu précédent, celui-ci est d'une
terrible difficulté, le groupe sous-jacent étant très
complexe.\footnote{Une bonne étude de ce groupe via le logiciel
  \textsc{Gap} est disponible à l'adresse
  \texttt{http://www.gap-system.org/Doc/Examples/rubik.html}. On y
  apprend entre autres que le groupe contient 43252003274489856000
  éléments.}

\subsubsection{Le taquin}
Plus facile, ce célèbre casse-tête du XIXème siècle présente un carré
de quatre cases par quatre, rempli par des petits carrés coulissants
portant les nombres de 1 à 15, la seizième case étant un espace vide,
qui permet d'y faire coulisser un des petits carrés adjacents. Partant
d'une configuration quelconque, le but du jeu est de rétablir la
configuration d'origine où les nombres sont alignés dans l'ordre, la
case vide restant en bas à droite.

Si ce jeu paraît très similaire au casse-tête du cube hongrois,
l'asymétrie causée par la présence d'une case vide fait en réalité que
sa meilleure formalisation ne s'exprime pas en termes de groupe mais
en termes de groupoïde, une variante plus faible de la notion de
groupe, qui a été présentée au \S\ref{ss:cat}.

Pour exprimer le taquin comme un jeu de factorisation sur un
groupoïde, il faut d'abord séparer l'ensemble $U$ des configurations
en 16 sous-ensembles $U_{i,j}$, suivant la position de la case vide
--- l'ensemble $U_{i,j}$ est l'ensemble des configurations pour
lesquelles la case vide se trouve en position $(i,j)$, $i$ et $j$
étant des entiers entre un et quatre représentant les coordonnées des
cases.  Dans le vocabulaire du \S\ref{ss:cat}, l'ensemble $X$ des
<<sommets du groupoïde>> est l'ensemble des 16 positions des cases du
jeu, et les <<flèches du groupoïde>> --- les opérations représentant
les manières de mélanger les petits carrés du jeu --- sont des flèches
allant d'une case de $X$ à une autre, ces deux cases étant les
positions respectives de la case vide avant et après application de
l'opération.  Dans cette description, la différence fondamentale avec
un groupe est que chaque opération ne peut s'appliquer qu'à certaines
configurations, celles dont la case vide se trouve au bon endroit.

Les mouvements fondamentaux (les générateurs) sont au
nombre\footnote{24 si on ne compte qu'une seule fois un mouvement et
  son inverse.} de 48 et s'expriment comme <<déplacer le petit carré en
position $(i_1,j_1)$ sur la case vide en position adjacente
$(i_2,j_2)$>>.  Comme la case vide se déplace alors de $(i_2,j_2)$ à
$(i_1,j_1)$, ces mouvements sont des flèches allant de la case
$(i_2,j_2)$ à la case $(i_1,j_1)$.

En plus d'utiliser un groupoïde au lieu d'un groupe, le jeu de taquin
présente une autre entorse à l'énoncé de jeu de factorisation tel
qu'il a été présenté au \S\ref{ss:factor}: à partir d'une
position de départ donnée, il est possible d'atteindre seulement la
moitié des configurations imaginables.\footnote{En réalité, le cube de
  Rubik ne permet pas non plus de réaliser toutes les configurations
  possibles, sans quoi il serait plus facile à résoudre.} Par
exemple, les mouvements appelés <<transpositions>>, qui permuteraient
deux des petits carrés en laissant tous les autres en place, sont
impossibles.  Une telle particularité ajoute à l'intérêt du jeu, et
les paragraphes suivants présenteront des idées de jeu pour exploiter
les situations de ce genre.

\subsection{Autres principes de jeu}\label{ss:varj}
Les trois ingrédients pour réaliser un jeu de théorie des groupes sont: 
\begin{itemize}
\item Un groupe (abstrait). Il n'est toutefois pas nécessaire de
  découvrir un nouveau groupe créer un jeu intéressant. Les groupes
  commutatifs sont en revanche à éviter, car ils ne permettent pas de
  construire un jeu basé sur l'ordre subtil d'opérations à effectuer.

  Les groupes possédant des identités plus complexes sont au contraire
  très intéressants car la découverte de cette identité révèle souvent
  une stratégie de gain. Par exemple le groupe de Heisenberg (voir le
  jeu~\ref{jeu:cocc}, page~\pageref{jeu:cocc}) est \emph{nilpotent};
  il satisfait l'identité $xyzyx=yxzxy$.

\item Un ensemble $U$ de configurations, faisant du groupe $G$ un
  groupe concret, pouvant être représenté de manière visuelle. Là
  aussi, les mathématiciens ont conçu beaucoup de représentations
  visuelles des groupes qu'ils étudient, car celles-ci leur permettent
  souvent d'améliorer la compréhension de la structure des groupes.

\item Le principe de jeu. Si la plupart des jeux existants basés sur
  la théorie des groupes sont de simples jeux de factorisation, il est
  pourtant possible de concevoir diverses variantes qui changent le
  caractère du jeu. Grâce aux grand nombres de groupes concrets
  intéressants qui pourront être utilisés pour chacun de ces jeux,
  l'invention de ces variantes donnera un réservoir quasi inépuisable
  de jeux.
\end{itemize}

Les sept descriptions qui suivent (et leurs variantes) sont des
formalisations de jeux existants, données plus ou moins dans un ordre
croissant de complexité. Le dernier paragraphe présentera au moins une
version concrète de chacun de ces jeux abstraits. Le contexte général
est toujours celui d'un groupe concret $G$ agissant transitivement sur
un ensemble $U$.

\begin{enumerate}
\item Un \emph{jeu de calcul mental} consiste en la donnée d'un
  élément $u_0\in U$ et d'un nombre fini d'éléments $g_1,\dots,g_n\in
  G$. Le but du jeu est de trouver l'élément $u=g_n\cdot g_{n-1}\cdots
  g_1\cdot u_0$.

  Dans l'exemple des nombres entiers donné au \S\ref{ss:factor}, le
  jeu de calcul mental consiste en la donnée des nombres 2, 2, 3 et 7,
  et le but du jeu est de calculer le produit de ces quatre nombres
  sans papier ni crayon.

\item Un \emph{jeu de calcul inverse} consiste en la donnée d'un
  élément $u_f\in U$ et d'un nombre fini d'éléments $g_1,\dots,g_n\in
  G$. Le but du jeu est de trouver l'élément $u\in U$ tel que
  $u_f=g_n\cdot g_{n-1}\cdots g_2\cdot g_1\cdot u$.

  Ce jeu demandant de raisonner <<à l'envers>> a été utilisé par
  Piaget pour tester à quel âge un enfant réalise les problèmes liés à
  la non-commutativité (voir \S\ref{ss:cat}) de l'enchaînement des
  opérations. Curieusement, \fcite{l'inversion des opérations\dots
    n'est comprise qu'à cet âge de 7-8 ans toujours retrouvé pour ce
    qui est de la découverte de la réversibilité.}{piaget:ar:inv}{52}
  En effet, la manière la plus facile (mais pas toujours la seule) de
  trouver la solution $u$ est d'inverser tous les éléments
  $g_1,\dots,g_n$ et de retrouver $u$ comme $g_1^{-1}\cdot
  g_2^{-1}\cdots g_n^{-1}\cdot u_f$, en prenant bien soin d'inverser
  l'ordre des opérations, alors que de jeunes enfants essaient souvent
  de calculer $g_n^{-1}\cdots g_2^{-1}\cdot g_1^{-1}\cdot u_f$, ce qui
  peut donner un résultat incorrect lorsque le groupe $G$ n'est pas
  commutatif.

  Pour les jeux qui suivent, la donnée supplémentaire est celle de $S$,
  un ensemble générateur du groupe $G$.
\item Un \emph{jeu de factorisation} consiste en la donnée de deux
  élément $u_0$ et $u_f \in U$. Le but du jeu est de choisir
  $s_1,s_2,\dots,s_n$ des éléments de $S$ de façon à ce que
  $u_f=s_n\cdot s_{n-1}\cdots s_2\cdot s_1\cdot u_0$. Les étapes
  successives $s_1\cdot u_0,s_2\cdot s_1\cdot u_0,\dots$ sont visibles
  par le joueur et la partie est gagnée dès que l'une de ces étapes
  est le but $u_f$.

  Ce jeu possède deux variantes:
  \begin{itemize}
  \item Une variante dite \emph{aveugle} qui le combine au \emph{jeu
      de calcul mental}, car les étapes intermédiaires ne sont pas
    visibles. Forcé de mentaliser ces étapes intermédiaires, le joueur
    sélectionne une séquence d'éléments de $S$, soumet sa séquence,
    puis observe le déroulement de ces opérations sans pouvoir
    intervenir.  La partie est gagnée si le but est atteint à la
    dernière étape exactement.

  \item Une variante dite \emph{contrainte} où chaque générateur ne
    peut être utilisé qu'un nombre limité de fois. Le joueur reçoit
    des cartes en nombre donné par une fonction $f:S\to\mathbb{N}$. Il
    doit exprimer $u_f$ comme $s_n\cdot s_{n-1}\cdots s_1\cdot u_0$ sans
    utiliser plus de cartes qu'il n'en a reçues, c'est-à-dire de façon
    à ce que pour tout $s\in S$ on ait $\#\{i|s_i=s\}\le f(s)$.
    Alternativement et plus simplement, la contrainte peut porter
    uniquement sur le nombre total de cartes qu'il est permis
    d'utiliser. Il est possible au joueur de reprendre son coup,
    c'est-à-dire d'annuler son dernier déplacement en reprenant en
    main la carte jouée. Bien sûr, la répartition des cartes doit être
    choisie de façon à ce qu'une solution soit possible.
  \end{itemize}

\item Un \emph{jeu de combinaison} est identique à un jeu de
  factorisation, sauf que l'élément but $u_f$ est inconnu du joueur,
  tel la combinaison d'un coffre-fort. Très différente, la stratégie
  sera alors de parvenir à parcourir tout l'ensemble $U$ en utilisant
  l'ensemble générateur $S$.

  Ici aussi, le jeu de combinaison possède une variante aveugle (à la
  différence près qu'ici la partie est gagnée si le but est atteint
  lors d'une étape intermédiaire) et une variante contrainte, qui ne
  sont pas mutuellement exclusives.

\item Un \emph{jeu de combinaison dévoilée} inclut la donnée
  supplémentaire de composantes d'une configuration. On suppose ici
  qu'on a fixé une famille $U_1,\dots,U_m$ de sous-ensembles de $U$. À
  chaque action du joueur, la configuration résultante $u$ est
  comparée à la configuration but $u_f$, et toutes les parties $U_i$
  qui contennent à la fois $u$ et $u_f$ sont révélés au joueur, et ces
  parties restent ensuite visibles jusqu'à la fin du jeu. Par exemple,
  si $U$ est l'ensemble des nombres à trois chiffres, on peut prendre
  $U_{1,0},\dots U_{3,9}$ les ensembles $U_{ij}$ des nombres ayant un
  $j$ à la position $i$. Ce jeu évoque le travail d'un spécialiste de
  l'ouverture des coffres-forts, qui est capable de repérer à
  l'oreille lorsqu'un des chiffres de la combinaison qu'il essaie de
  retrouver est à sa place.

  Dans la variante dite \emph{memory}, les composantes révélées au
  joueur ne sont plus visibles dès que la configuration courante a
  changé; le joueur doit donc mémoriser ces composantes révélées.

  Le \emph{jeu de combinaison dévoilée} possède également une variante
  contrainte et une variante aveugle.

  Pour introduire les deux derniers jeux de ce paragraphe, il est
  nécessaire d'introduire une notion très importante de théorie des
  groupes, celle de sous-groupe. Un \emph{sous-groupe} $H$ d'un groupe
  (concret ou abstrait) $G$ est un sous-ensemble de $G$ qui est
  lui-même un groupe pour la restriction de la multiplication de $G$
  (voir \S\ref{ss:cat}), ce qui revient à dire qu'il
  satisfait les trois propriétés suivantes:
  \begin{enumerate}
  \item Le produit de deux éléments de $H$ est encore dans $H$, i.e.\
    $h_1,h_2\in H\Longrightarrow h_1h_2\in H$;
  \item L'inverse d'un élément de $H$ est encore dans $H$, i.e.\ $h\in
    H\Longrightarrow h^{-1}\in H$;
  \item L'élément neutre de $G$ est dans $H$, i.e.\ $1\in H$.
  \end{enumerate}
  
  Les deux derniers jeux utilisent les mêmes données que les jeux qui
  précèdent, mais avec en plus un sous-groupe $H$ du groupe concret
  $G$. Le sous-groupe $H$ agit également sur l'ensemble $U$ (puisque
  les éléments de $H$ sont des éléments de $G$), mais $H$ contenant
  moins d'éléments que $G$, il ne sera pas toujours possible de passer
  d'un élément de $U$ à un autre en utilisant les opérations de $H$,
  c'est-à-dire que l'action de $H$ sur $U$ sera \emph{non transitive}
  (voir \S\ref{ss:factor}). On demande tout de même que
  le nombre d'orbites de $H$ soit fini.\footnote{C'est par exemple
    vérifié si le sous-groupe $H$ est d'\emph{indice fini} dans G.}

\item Un \emph{jeu de factorisation impossible} est la traduction d'un
  \emph{jeu de factorisation} pour cette action non transitive, la
  seule différence étant que l'ensemble $S$ est un ensemble générateur
  de $H$ et non un ensemble générateur de $G$. Les configurations de
  départ ($u_0$) et de but ($u_f$) sont toujours choisies au hasard,
  ce qui signifie que le jeu possèdera parfois une solution (lorsque
  $u_0$ et $u_f$ se trouvent dans la même orbite) et parfois aucune
  (lorsque $u_0$ et $u_f$ se trouvent dans des orbites différentes).
  Dans ce dernier cas, le joueur gagne uniquement en déclarant
  <<mission impossible>>.  Bien entendu, le joueur perd s'il déclare
  <<mission impossible>> alors que le jeu possédait une solution.

  D'après les recherches de Piaget, la prise en compte d'un élément
  d'impossibilité logique ne devrait être maîtrisable que par des
  enfants plus âgés (environ 10 ans). Avant, \fcite{les sujet observés
    commencent par essayer\dots et si la moitié des sujets finit par
    se rallier à l'impossibilité, les autres admettent qu'ils n'y
    parviendront pas, mais pensent qu'un adulte astucieux y arriverait
    sans doute.}{piaget:rc:jeu}{49}

\item Un \emph{jeu de factorisation à choix} consiste en la donnée
  d'un élément $u_0\in U$ et d'un sous-ensemble $T\subset U$ qui est
  une \emph{transversale} pour l'action de $H$ sur $U$, c'est-à-dire
  qui contient exactement un élément de chaque orbite --- cet ensemble
  $T$ est fini car nous avons supposé le nombre d'orbites fini. Le but
  du jeu est de trouver $s_1,s_2,\dots,s_n\in S$ tels que $s_n\cdot
  s_{n-1}\cdots s_2\cdot s_1\cdot u_0$ soit un élément de la
  transversale $T$.  Autrement dit, le jeu possède en apparence
  plusieurs buts (tous les éléments de $T$ en sont), mais dont en
  réalité un seul est atteignable: l'élément de $T$ qui est dans la même
  orbite que $u_0$.
\end{enumerate}

Ces deux derniers jeux sont peut-être les plus intéressants, en ce
qu'ils amènent à une compréhension intuitive de la notion de
sous-groupe. Ils peuvent également se jouer dans la variante aveugle,
mais il semble par contre pédagogiquement maladroit d'en faire des
variantes contraintes, car cela conduirait trop facilement le joueur à
confondre une impossibilité <<fortuite>> d'atteindre une solution par
manque de cartes, et l'impossibilité <<structurelle>> de l'atteindre
car elle se trouve dans une orbite différente.

\subsection{Quelques jeux, implémentables ou déjà implémentés sur
  ordinateur}\label{ss:implementes}
L' utilisation de l'ordinateur est un apport très important pour
réaliser les jeux décrits au \S\ref{ss:varj}. Premièrement,
elle permet d'accéder à une plus grande variété de représentations
visuelles des groupes concrets, se libérant des contraintes liées à la
réalisation technique des jeux. Par exemple, il est à peu près certain
que la création du cube de Rubik a demandé à son auteur beaucoup plus
d'ingéniosité technique que d'imagination mathématique. D'une façon
générale, dans un casse-tête physique du type <<factorisation>>, le défi
du constructeur est de restreindre le joueur à un ensemble générateur,
c'est-à-dire un ensemble de mouvements possibles, qui soit
suffisamment non trivial pour rendre le jeu intéressant; avec
l'informatique par contre, rien de plus facile que de ne définir comme
commandes que certains boutons qui auront l'action désirée sur la
configuration représentée à l'écran.

Deuxièmement, seule l'utilisation d'un ordinateur permet facilement de
cacher au joueur certaines informations, rendant ainsi possible les
jeux de combinaison et les variantes aveugles ou de calcul mental.
Finalement, un jeu sur ordinateur ne permet pas de tricher lorsque le
jeu n'est pas uniquement défini par ses contraintes physiques, et que
d'autres règles doivent être énoncées. Ceci est par exemple le cas
pour la variante de factorisation contrainte, qui fait intervenir des
<<cartes>>.

La réalisation de ces jeux présente deux bénéfices importants: d'une
part, elle donne des exemples que la théorie des groupes peut surgir
naturellement à partir de situations très simples. Mais
principalement, elle permet de proposer à des enfants un cadre
purement ludique, dans lequel ils seront conduits à percevoir
intuitivement la présence d'un groupe comme une entité cachée, et à
s'entraîner à l'utilisation d'un raisonnement opératoire.

Regroupés en trois grandes familles, voici quelques exemples de jeux
qui ont été pensés pour être implémentés sur ordinateur. Les
jeux~\ref{jeu:pl}, \ref{jeu:ps}, \ref{jeu:coffre} et~\ref{jeu:elef}
ont déjà été programmés, et tous les autres ne demandent qu'à l'être.
Les principes de jeu ayant déjà été exposés au \S\ref{ss:varj}, la
description de ces jeux mettra principalement l'accent sur la
description du groupe concret --- les configurations visibles et leurs
opérations sous-jacentes --- ainsi que sur l'interface proposée au
joueur.

La première famille de jeux utilise les groupes de permutations de $n$
objets, groupes qui sont parmi les plus fondamentaux utilisés en
mathématiques. Le groupe des permutations de $n$ objets est le groupe
concret de toutes les permutations possibles de $n$ objets différents
occupant $n$ positions distinctes. Dans notre cas, les $n$ objets
différents seront des billes de couleurs différentes\footnote{Bien
  sûr, si cela est jugé plus attractif, les billes de couleurs peuvent
  être remplacées par des motifs plus plaisants.}, qui seront souvent,
mais pas forcément, disposées en ligne. Avec par exemple une bille
jaune J, une bille rouge R et une bille bleue B, il y a six
configurations possibles, qui sont JRB, JBR, RJB, RBJ, BJR et BRJ. Les
opérations agissant sur de telles configurations sont décrites par
exemple comme <<inverser les positions de la deuxième et de la
troisième bille>>, cette opération transformant JRB en JBR.

\begin{enumerate}
\item Le \emph{jeu de la machine infernale} est un \emph{jeu de calcul
    mental} sur trois (ou plus) billes de couleur. Les trois billes se
  présentent dans un certain ordre à l'entrée d'un tube, puis leurs
  couleurs deviennent invisibles dès qu'elles sont rentrées dans le
  tube --- on voit alors des billes uniformément grises. A plusieurs
  endroits du tube, des dispositifs en forme de roue, et pouvant
  accueillir deux ou trois billes, permettent d'inverser l'ordre de
  celles-ci en effectuant une moitié de tour complet, ceci étant
  observé par le joueur. Juste avant que les billes ressortent du tube
  par une autre ou par la même extrémité, le joueur doit prédire dans
  quel ordre les billes vont ressortir.

  Dans une expérience de Piaget, décrite au début du
  chapitre~\ref{s:math}, ce jeu a été utilisé dans sa forme la plus
  simple, où les trois billes subissent exactement une inversion de
  leur ordre, et où le joueur doit se contenter de prédire quelle est
  la première bille qui sortira du tube.

\item\label{jeu:pl} Le \emph{jeu des permutations linéaires} est un
  \emph{jeu de factorisation} sur cinq (ou plus) billes de couleurs
  disposées en ligne, avec comme déplacements possibles (système
  générateur) les quatre (ou plus) opérations les plus élémentaires, à
  savoir celles qui permutent les positions de deux billes contiguës.
  Comme dans le jeu précédent, ces opérations peuvent très bien être
  représentées par la rotation d'une roue. En utilisant ces
  opérations, le joueur doit recomposer une configuration qui lui est
  donnée comme but.

\item\label{jeu:ps} Le \emph{jeu des permutations cycliques} est un
  \emph{jeu de factorisation} sur cinq (ou plus) billes de couleur
  disposées sur un plateau circulaire --- et formant donc un pentagone
  régulier --- avec cette fois seulement deux générateurs: le premier
  générateur permute les positions des deux billes situées au sommet
  du cercle, tandis que le deuxième --- pouvant être représenté par
  une manivelle faisant tourner le plateau circulaire --- fait tourner
  le cercle entier d'un cinquième de tour dans le sens désiré, chaque
  bille prenant alors la position de sa voisine.

\item Pour les joueurs plus avancés, le \emph{jeu des permutations
    paires} est un jeu de factorisation impossible sur sept billes de
  couleur, disposées sur deux plateaux circulaires se chevauchant; le
  premier plateau porte trois billes (formant un triangle équilatéral)
  et le second cinq billes (formant un pentagone régulier), une des
  sept billes se trouvant à l'intersection des deux plateaux et
  faisant à la fois partie du triangle et du pentagone. Le joueur doit
  atteindre la configuration but en se servant de deux manivelles
  permettant de faire tourner les plateaux. La difficulté et l'intérêt
  supplémentaires de ce jeu est que ces deux opérations ne sont pas
  des générateurs du groupe de permutations, mais ne permettent
  d'atteindre que l'exacte moitié des configurations possibles,
  correspondant au sous-groupe des permutations dites \emph{paires}.
  Si les configurations de départ et de but sont choisies au hasard,
  le joueur devra donc une fois sur deux déclarer <<mission
  impossible>>.

\item\label{jeu:coffre} Le \emph{jeu du coffre} est un \emph{jeu de
    combinaison} sur quatre billes de couleurs disposées en ligne. Les
  générateurs sont les mêmes que dans le \emph{jeu des permutations
    linéaires}, mais la combinaison à réaliser est inconnue et sert à
  ouvrir un coffre au trésor. Pour faire comprendre le but au joueur,
  les billes font partie d'un dispositif scellé sur un coffre, où des
  petites roues similaires à celle du \emph{jeu de la machine
    infernale} permettent de permuter deux billes contiguës.

  La deuxième famille de jeux est basée sur un groupe plus familier
  aux non-mathématiciens, celui des déplacements dans le plan. Le
  plateau de jeu y sera toujours la grille quadrillée décrite dans le
  premier exemple du \S\ref{ss:class}, sur laquelle un pion
  représenté par un dessin d'animal pourra être déplacé.

\item Le \emph{jeu des moutons} est un \emph{jeu de factorisation
    contrainte et aveugle}. Sur différentes cases du plateau de jeu,
  on trouve un mouton vu d'en haut, dont la tête est orientée vers une
  des quatre directions possibles, une bergerie dont l'entrée se fait
  par l'une des quatre directions possibles, ainsi que des obstacles
  représentant des cases infranchissables. Le joueur reçoit quelques
  cartes, sur chacune desquelles figure un des trois déplacements
  suivants du mouton: <<avancer>>, <<pivoter sur sa gauche>> ou
  <<pivoter sur sa droite>>. Le but du jeu est de composer une
  séquence de cartes qui dirigeront le mouton vers sa bergerie en
  évitant les obstacles.  Lorsque le joueur a choisi sa séquence de
  cartes, il soumet sa solution, qui est appliquée par le programme;
  il gagne si le mouton entre dans sa bergerie par le bon côté au
  cours de la solution.

\item Le \emph{jeu des moutons programmés} est un \emph{jeu de calcul
    inverse} utilisant le même cadre que le jeu précédent. Sont ici
  donnés la bergerie et une séquence de flèches décrivant les
  mouvements qui seront effectués par le mouton (toujours parmi
  <<avancer>>, <<gauche>> et <<droite>>.). Placer des obstacles n'est
  ici pas nécessaire. Le joueur choisit un des quatre moutons qui ont
  été représentés à l'écran dans chacune des quatre orientations
  possibles, puis clique sur une case, qui sera la case de départ de
  son mouton. Celui-ci effectue alors sa séquence de programmation, et
  la partie est gagnée uniquement si le mouton entre dans sa bergerie
  avec la dernière flèche de sa séquence --- qui sera toujours une
  flèche <<avancer>>.

\item\label{jeu:cocc} Le \emph{jeu de la coccinelle} fait intervenir
  un groupe plus complexe, appelé \emph{groupe de Heisenberg}. La
  coccinelle est au départ de couleur unie, mais ses déplacements
  verticaux vont lui ajouter ou retirer des points sur le dos. Les
  quatre déplacements possibles (générateurs) sont simplement Sud,
  Ouest, Nord et Est, comme dans le premier jeu du \S\ref{ss:class}.
  Le plateau de jeu est de sept cases sur sept, et chaque colonne sauf
  la première est repérée par un dé à jouer indiquant de 1 pour la
  deuxième colonne à 6 pour la septième.  La coccinelle part de la
  case en bas à gauche (à l'extrême sud-ouest). Un déplacement en
  direction de l'est ou de l'ouest laisse la coccinelle inchangée; un
  déplacement en direction du nord ajoute à la coccinelle autant de
  points que l'indique la colonne sur laquelle elle se trouve; un
  déplacement en direction du sud retire à la coccinelle autant de
  points que l'indique la colonne sur laquelle elle se trouve; si la
  coccinelle ne porte pas assez de points pour cela, le déplacement
  vers le sud est impossible. Le jeu est une \emph{factorisation
    contrainte}, où le nombre total de déplacements de la coccinelle
  est limité à un maximum de 16; le but est de ramener la coccinelle
  sur sa case de départ en ayant accumulé le plus de points possible.
  Le nombre de points maximal est de 16 et s'obtient en faisant quatre
  pas vers l'est, quatre pas vers le nord, quatre pas vers l'ouest
  puis quatre pas vers le sud; revenir avec 16 points est donc le but
  du jeu, mais on peut bien sûr mesurer le score du joueur par le
  nombre de points obtenu.

  Les deux derniers jeux de la famille des déplacements dans le plan
  sont des \emph{jeux de factorisation à choix}, utilisant des groupes
  connnus des mathématiciens sous la dénomination générale de
  \emph{groupes cristallographiques}.

\item\label{jeu:elef} Dans le \emph{jeu des éléphants réfléchis}, joué
  sur un quadrillage de six cases par six\footnote{Ou plus grand, mais
    le nombre de cases doit être pair.}, la figurine représente un
  éléphant vu de profil, c'est-à-dire ne présentant aucun axe de
  symétrie. Dans sa position de départ, l'éléphant est représenté à
  l'endroit, avec sa trompe vers la droite. Ici encore, les
  déplacements possibles sont Sud, Ouest, Nord et Est. Toutefois,
  chaque déplacement de l'éléphant a un <<effet secondaire>>: l'image
  de l'éléphant est retournée, exactement comme si toutes les arêtes
  séparant les cases étaient des miroirs. Si l'éléphant se déplace
  vers l'ouest ou vers l'est, il se retrouvera donc avec la trompe à
  gauche. Si maintenant il se déplace vers le nord ou vers le sud, il
  se retrouvera sur le dos. Avec ces effets secondaires, le groupe
  engendré par l'ensemble générateur \{Sud, Ouest, Nord, Est\} est
  alors plus petit d'un facteur quatre que le groupe de tous les
  mouvements envisageables, celui qui permettrait d'obtenir n'importe
  quelle position de l'éléphant. Ceci signifie que l'ensemble des
  configurations possibles est formé de quatre orbites, et permet de
  faire du jeu des éléphants réfléchis un \emph{jeu de factorisation à
    choix}: le but du jeu est d'amener l'éléphant <<à l'endroit>>,
  c'est-à-dire sur ses pattes et regardant vers la droite, sur l'une
  parmi quatre cases-but disposées aux quatre coins du quadrillage ---
  ces quatre cases se trouvant dans les quatre orbites différentes.

  Optionnellement, il est possible de fournir une aide au joueur en
  coloriant la grille en quatre couleurs, suivant la parité du numéro
  de ligne et du numéro de colonne. Ce coloriage représentant
  exactement la stratification en orbites des cases du quadrillage, il
  permet de déterminer plus facilement quelle case-but doit être
  visée: si par exemple l'éléphant se trouve à l'endroit sur une case
  jaune, cela signifie qu'il doit viser la case d'arrivée jaune.
  D'autre part, une variante peut-être plus amusante propose quatre
  éléphants au lieu d'un, qui sont disposés au départ de façon à être
  rangés chacun dans une des quatre cases-but.

\item En apparence quasiment identique, le \emph{jeu des éléphants
    tournants} introduit une composante supplémentaire. La
  présentation en est similaire au \emph{jeu des éléphants réfléchis},
  mais avec un petit cercle sur chaque intersection des lignes,
  représentant un centre de rotation. Les seuls déplacements possibles
  consistent à faire tourner d'un quart de tour les quatres cases se
  trouvant autour de l'intersection choisie, cette rotation ayant bien
  sûr effet sur l'orientation de l'éléphant. L'éléphant a maintenant
  huit orientations possibles, et pourtant, tout comme dans le
  \emph{jeu des éléphants réfléchis}, l'ensemble des configurations
  possibles est toujours formé de seulement quatre orbites,
  correspondant au coloriage décrit plus haut. Ceci s'explique par la
  présence d'un mouvement d'aller-retour qui permet de renverser
  complètement l'éléphant: faire d'abord tourner d'un quart de tour
  dans le sens des aiguilles d'une montre le carré de quatre cases
  dont l'éléphant se trouve sur la case nord-ouest (l'éléphant se
  déplace d'une case vers l'est), puis faire tourner, toujours dans le
  sens des aiguilles d'une montre, le carré de quatre cases dont
  l'éléphant se trouve maintenant sur la case sud-est. Ce deuxième
  mouvement ramène l'éléphant sur sa case originale, mais comme il a
  tourné deux fois dans le sens des aiguilles d'une montre, il se
  retrouve renversé à $180^o$. La découverte et l'identification de
  tels mouvements de renversement est essentielle pour réussir ce jeu
  de façon sûre.\footnote{Pour les mathématiciens, ce jeu illustre un
    résultat élémentaire: le nombre d'orbites (4) multiplié par le
    cardinal d'un stabilisateur (2) est égal à l'indice du sous-groupe
    (8).}

  La variante avec quatre éléphants au lieu d'un présente l'attrait
  supplémentaire qu'il sera possible de faire tourner plusieurs
  éléphants en même temps, si ceux-ci se trouvent sur le même carré de
  quatre cases.

  La troisième famille de jeu utilise des puzzles constitués de pièces
  simples, à placer de façon à obtenir un motif dont le modèle est
  visible. Tous ces puzzles auront la particularité que toutes leurs
  pièces occupent en permanence toutes les cases possibles, mais sont
  disposées au départ dans un mauvais ordre. Un mouvement possible
  représentera donc toujours le mouvement simultané de plusieurs
  pièces, exactement comme un jeu de taquin sans case vide, l'absence
  de case vide permettant de faire intervenir de vrais groupes au lieu
  de groupoïdes (voir \S\ref{ss:class}).  Pour ne pas
  faire interférer des problèmes de théorie des groupes et des
  problèmes de reconnaissance visuelle, on choisira généralement un
  motif facile à reconstituer, voire éventuellement de simples nombres
  à remettre dans l'ordre, comme dans le jeu de taquin.

\item Le \emph{puzzle carré libre} est un \emph{jeu de factorisation}
  sur un carré de trois cases sur trois (ou plus), avec une grande
  liberté de mouvement: en cliquant successivement sur deux cases, il
  est possible de faire permuter n'importe quelle paire de pièces en
  laissant toutes les autres en place. Ces nombreux générateurs sont
  ceux qui permettent de résoudre le jeu le plus facilement.

\item Le \emph{puzzle carré taquin} limite les mouvement du jeu
  précédent aux permutations des paires de pièces adjacentes,
  effectuées en cliquant entre les deux cases. Il s'agit d'une version
  simplifiée du jeu du taquin (voir \S\ref{ss:class}), où
  il n'est pas nécessaire de s'inquiéter de la position de la case
  libre.

\item Le \emph{puzzle carré avec rotations} est toujours un jeu de
  factorisation sur le même carré de quatre cases sur quatre, dont les
  générateurs sont les rotations des carrés de quatre cases. Il s'agit
  des mêmes mouvements que ceux décrits dans le \emph{jeu des
    éléphants tournants}, mais sans faire tourner les pièces sur
  elles-même. Il est possible d'atteindre ainsi n'importe quelle
  configuration, mais c'est nettement plus difficile que dans les deux
  premiers jeux.

\item Le \emph{puzzle carré cristallographique} utilise les mêmes
  rotations que le jeu précédent, mais cette fois ces rotations font
  également tourner les les pièces sur elles-même, exactement comme
  dans le \emph{jeu des éléphants tournants}. Cette fois il n'est pas
  possible d'atteindre n'importe quelle configuration, le nombre
  d'orbites étant de quatre. Selon les goûts, il s'agira donc soit
  d'un \emph{jeu de factorisation impossible}, soit d'un \emph{jeu de
    factorisation à choix}, avec quatre modèles dont un seul est
  réalisable.

\item Le \emph{puzzle hexagonal} est composé de sept hexagones
  réguliers, un hexagone central et six hexagones disposés tout autour
  comme dans un nid d'abeille. Tout comme pour le \emph{puzzle carré
    avec rotations}, les mouvements possibles sont les rotations
  autour des six intersections (les sommets de l'hexagone central),
  qui font ici tourner trois hexagones d'un tiers de tour, sans que
  ceux-ci tournent sur eux-mêmes. Exactement comme dans le \emph{jeu
    des permutations paires}, ces mouvements ne permettent d'atteindre
  que la moitié des configurations, il s'agit donc d'un \emph{jeu de
    factorisation impossible}.

\item Le \emph{puzzle hexagonal corrigé} est un \emph{jeu de
    factorisation} identique au jeu précédent, sauf que deux des
  pièces hexagonales sont identiques. Cette petite modification est en
  effet suffisante pour que la factorisation ne soit plus jamais
  impossible.\footnote{Cette variante est directement inspirée d'un
    ingénieux jeu de taquin utilisant, au lieu des nombres de un à
    quinze, les lettres de "DEAD PIGS WONT FLY" réparties ainsi sur
    les quatre lignes du jeu. Le but est alors de transformer le
    message en "DEAD PIGS FLY TOWN". L'astuce est qu'ici le jeu admet
    une solution, en permutant les deux D du mot DEAD.}

\item Le \emph{puzzle hexagonal cristallographique} est la version du
  \emph{puzzle hexagonal} où les hexagones tournent également sur
  eux-mêmes. Les orbites étant ici au nombre de six, il semble
  judicieux d'en faire un \emph{jeu de factorisation à choix} en
  proposant six modèles dont un seul est réalisable.

\item Le \emph{puzzle triangulaire} est composé de huit triangles
  équilatéraux formant un losange. Six des triangles sont disposés
  autour d'un point central pour former un hexagone régulier, auquel
  on ajoute deux triangles sur des côtés opposés de l'hexagone
  central. Les mouvements possibles sont ici la rotation de l'hexagone
  (les triangles tournant également sur eux-mêmes), ainsi que les deux
  réflexions-miroirs qui permutent l'un des deux triangles extérieurs
  avec le triangle intérieur qui lui est adjacent.

  Il est facile de voir que si un triangle est ramené à sa position
  initiale par une suite de telles opérations, alors il reviendra
  forcément dans la même orientation; par ailleurs, les huit objets
  (les triangles) peuvent être permutés arbitrairement, et le groupe
  associé à ce jeu est le groupe symétrique $S_8$, à 40320 éléments.
  On peut considérer ce jeu comme un \emph{jeu de factorisation}, ou
  de \emph{factorisation impossible}, ou même d'\emph{isomorphisme}
  (voir ci-dessous).
\end{enumerate}

\subsection{Quelques pistes supplémentaires de jeux}
En plus des jeux présentés au \S\ref{ss:implementes}, il est possible
d'imaginer d'autres pistes à exploiter; il s'agira ici de jeux
incomplets dont le principe ou la présentation précise restera
partiellement à déterminer.

Les jeux nommés \emph{puzzles} dans le \S\ref{ss:implementes}
inspirent un prolongement qui ne surprendra guère les mathématiciens,
celui des \emph{pavages euclidiens}. A la différence de ces
\emph{puzzles}, un \emph{pavage euclidien} doit être vu comme une
division \emph{infinie} du plan entier en figures géométriques. Les
trois exemples qui ont été présentés sont le pavage par des carrés,
appelé quadrillage, le pavage par des triangles équilatéraux et la
pavage par des hexagones réguliers; il se trouve que ce sont les trois
seules manières de remplir le plan avec des figures toutes identiques,
ceci à déformation des pièces près. Il est en effet possible de
déformer les pièces de ces pavages jusqu'à les rendre méconnaissables,
comme l'a souvent fait l'artiste M.C. Escher dans ses peintures et
gravures. Ce qui reste par contre invariant par déformation des pièces
est le \emph{groupe d'isométries} sous-jacent du pavage.  Le
\emph{groupe d'isométries} d'un pavage est un sous-groupe de
l'ensemble des isométries du plan, c'est-à-dire l'ensemble de toutes
les transformations du plan qui laissent toutes les distances entre
les points inchangées\footnote{Ainsi, par conséquent, que tous les
  angles.}.  Le groupe des isométries du plan contient toutes les
translations, toutes les rotations, toutes les symétries axiales,
ainsi que les compositions de ces transformations. Les isométries du
pavage sont celles qui laissent celui-ci inchangé; par exemple, pour
un quadrillage, la rotation de 90 degrés autour d'un sommet déplace le
quadrillage sur lui-même, et il s'agit donc d'une isométrie de ce
pavage.  Les possibilités de pavage du plan sont innombrables
lorsqu'on permet d'utiliser plus d'une sorte de figure géométrique. À
chacun de ces pavages est donc associé un groupe que l'on appelle
\emph{groupe cristallographique}, qui permet en particulier de
reconnaître si deux pavages sont <<équivalents>>, c'est-à-dire que l'un
ne serait qu'une déformation de l'autre.  On peut imaginer des jeux
montrant l'un de ces pavages dont chaque sorte de figure serait d'une
couleur différente. Pour donner au joueur l'impression que le pavage
est réellement infini, il faudrait sans doute lui donner la possiblité
de naviguer dans celui-ci en le faisant défiler sur l'écran. Le jeu
pourrait alors être de manipuler ce pavage à l'aide de certaines
translations, rotations et symétries pour l'amener sur un modèle
donné.

Un autre piste est celle des groupes de déplacement non-euclidiens. À
deux dimensions, les <<géométries modèles>> sont au nombre de trois:
euclidienne, sphérique et hyperbolique. À trois dimensions, il y en a
huit; voir~\cite{thurston:3d} pour un exposé de ce sujet. Ceci
signifie qu'il y a huit espaces qui ressemblent localement à l'espace
usuel, mais où les notions de distance, et la structure globale sont
différentes. L'une de ces géométries, appelée \emph{Nil}, a déjà été
présenté de manière déguisée dans le \emph{jeu de la coccinelle} (voir
le \S\ref{ss:implementes}). L'espace sous-jacent peut y être identifié
à $\R^3$, mais les droites euclidiennes ne <<sont plus droites>> dans
cette nouvelle géométrie, et vice versa. Malgré les apparences, la
coccinnelle évoluait véritablement dans un espace à trois dimensions,
la troisième dimension étant repérée par le nombre de points qu'elle
portait sur son dos. Les règles régissant l'augmentation ou la
diminution de ce nombre de points sont en fait exactement les règles
qui déterminent la géométrie \emph{Nil}.  De façon similaire aux
pavages présentés plus haut, chacune de ces huit géométries de
l'espace est univoquement déterminée par son \emph{groupe
  d'isométries}, c'est-à-dire par le groupe des transformations de
l'espace qui préservent toutes les distances entre les points. Le
groupe d'isométries de l'espace euclidien standard est bien entendu le
groupe des déplacements introduit par Piaget et décrit au
\S\ref{ss:gd}.  La particularité des jeux de déplacements dans des
géométries non standard serait que, tout comme dans le cas de la
coccinnelle, certains déplacements seraient plus <<efficaces>> que
d'autres, et que le problème d'aller d'un point à un autre de l'espace
le plus rapidement possible ne consisterait plus à décrire une ligne
droite?  La solution à ce problème est, par définition, de parcourir
une \emph{géodésique} de la géométrie concernée. La manière la plus
ludique d'exposer cette problématique est sans doute de proposer un
jeu d'action en trois dimensions. Bien qu'on puisse arguer que de tels
jeux fassent surtout appel à des notions de géométrie, la présence
sous-jacente du groupe d'isométries y est aussi importante qu'elle l'a
été soulignée par Piaget dans ses travaux sur les groupes de
déplacements.

Au cours de cet article tout comme dans le discours de Piaget,
l'accent a souvent été mis sur la distinction entre une situation
concrète, présentée soit par une situation de la vie réelle, soit par
un jeu inventé de toutes pièces, et le groupe sous-jacent des
opérations applicables à cette situation ou à ce jeu, qui pour sa part
ne peut pas être appréhendé visuellement et dont la conscience demande
un effort d'abstraction. Un corollaire de cette dualité est qu'il est
tout à fait possible de rencontrer deux situations en apparence
totalement différentes, mais dont les groupes sous-jacents sont
identiques. En réalité, ces deux groupes ne sont pas identiques en
tant que groupes concrets, puisqu'ils agissent sur des ensembles
différents, mais ils ont la même structure, et par conséquent sont
identiques en tant que groupe abstraits. Le terme mathématique exact
est que ces deux groupes sont \emph{isomorphes}\footnote{Plus
  précisément, ceci signifie que les éléments des deux groupes peuvent
  être mis en correspondance bijective, et que cette correspondance
  préserve l'opération du groupe.}.  Lorsque deux situations concrètes
présentent le même groupe sous-jacent, cela signifie en pratique que
toute stratégie qui s'applique à la première situation s'applique
également à la deuxième, et vice versa. La capacité à identifier cet
\emph{isomorphisme de groupes} et à pouvoir transposer une stratégie
d'une situation à une autre est une composante extrêmement importante
de l'intelligence mathématique. Construire des jeux se basant sur la
reconnaissance de ces isomorphismes serait donc particulièrement
intéressant.  Un exemple très---peut-être trop---simple
d'\emph{isomorphisme de groupes} est donné par les deux situations
suivantes. Dans la première, un petit dessin carré ne présentant
aucune symétrie peut être tourné de 90 degrés dans le sens désiré, ou
alors réfléchi selon un axe vertical. Il s'agit de la représentation
standard du \emph{groupe des isométries du carré}. La deuxième
situation présente deux ampoules absolument identiques. Un
interrupteur permet d'allumer ou d'éteindre l'ampoule de gauche, et un
deuxième bouton permet de permuter les deux ampoules (par exemple, en
faisant tourner un plateau circulaire). Bien qu'il n'y paraisse pas,
le groupe sous-jacent à l'ensemble de configurations des deux ampoules
est isomorphe au groupe des isométries du carré.  On peut imaginer
deux types de jeux se basant sur deux ensembles de configurations dont
les groupes sous-jacents sont isomorphes. La présentation la plus
évidente divise l'écran d'ordinateur verticalement en deux, et
présente les deux configurations de part et d'autre. La manipulation
d'une des deux configurations provoque automatiquement la
transformation équivalente de l'autre configuration. Le défaut de cet
isomorphisme très explicite est que le joueur aura tendance à
manipuler la configuration avec laquelle il est le plus à son aise et
à négliger l'autre.  Une présentation beaucoup plus implicite serait
de construire deux fois le même jeu, de préférence relativement
difficile, avec ces deux ensembles de configurations différents.
Exposé au deuxième de ces jeux, le sujet pourrait alors profiter de
l'expérience accumulée avec le premier, mais il serait bien sûr tout à
fait possible qu'il ignore celle-ci et recommence son analyse à zéro.
Peut-être le mécanisme de jeu idéal pour mettre en évidence le
principe d'isomorphisme se situe-t-il quelque part entre la
présentation explicite et la présentation implicite.

\end{document}